\documentclass{article}

\usepackage{hyperref}

\usepackage{geometry}
\usepackage{graphicx}
\graphicspath{{figures/}}
\DeclareGraphicsExtensions{.png,.pdf}

\usepackage{mathtools}
\usepackage{amsmath,amsfonts,amsthm,amssymb}

\usepackage{tikz,pgfplots,pgfplotstable}
\usetikzlibrary{positioning,calc}
\usetikzlibrary{decorations.pathreplacing,calligraphy}
\usetikzlibrary{shapes.geometric}

\newtheorem{definition}{Definition}
\newtheorem{remark}{Remark}

\usepackage{multirow}
\usepackage{booktabs}
\usepackage{caption,subcaption}
\usepackage{cite}
\usepackage{nicefrac}

\newcommand{\NumFlux}{\mathcal{F}}
\newcommand{\dx}{\mathrm{d}x}
\renewcommand{\d}{\mathrm{d}}

\newcommand{\R}{\mathbb{R}}

\newcommand{\Popt}{P_{\text{\sf opt}}}
\newcommand{\Poly}[1]{\mathbb{P}^{#1}}
\newcommand{\Sopt}{\mathcal{S}_{\text{\sf opt}}}

\newcommand{\ca}[1]{\overline{#1}}
\newcommand{\DT}{\mathrm{\Delta} t}

\newcommand{\CW}{\mathsf{CWZb}3}
\newcommand{\Q}{\mathsf{Q}3}
\newcommand{\QI}[1]{\mathsf{Q}3_{\scriptscriptstyle{\mathcal{I}#1}}}

\newcommand{\QP}{\mathsf{Q}3\mathsf{P}1}
\newcommand{\QPMOOD}{\mathsf{Q}3\mathsf{P}1_{\mathsf{MOOD}}}

\newcommand{\TotDer}[2]{\frac{\mathrm{d}{#1}}{\mathrm{d}{#2}}}

\newcommand{\bigvec}[1]{\mathfrak{#1}}

\title{Quinpi: Integrating stiff hyperbolic systems with implicit high order finite volume schemes}

\author{
	Gabriella Puppo\footnote{Dipartimento di Matematica -- Sapienza, Universit\`{a} di Roma; P.le Aldo Moro, 5 -- 00185 Roma (Italy); gabriella.puppo@uniroma1.it}
	\quad Matteo Semplice\footnote{Dipartimento di Scienza e Alta Tecnologia -- Universit\`{a} dell'Insubria; Via Valleggio, 11 -- 22100 Como (Italy); matteo.semplice@uninsubria.it}
	\quad Giuseppe Visconti\footnote{Dipartimento di Matematica -- Sapienza, Universit\`{a} di Roma; P.le Aldo Moro, 5 -- 00185 Roma (Italy); giuseppe.visconti@uniroma1.it}
	}
\begin{document}

\maketitle

\begin{abstract}
	Many interesting physical problems described by systems of hyperbolic conservation laws are stiff, and thus impose a very small time-step because of the restrictive CFL stability condition. In this case, one can exploit the superior stability properties of implicit time integration which allows to choose the time-step only from accuracy requirements, and thus avoid the use of small time-steps. We discuss an efficient framework to devise high order implicit schemes for stiff hyperbolic systems without tailoring it to a specific problem. The nonlinearity of high order schemes, due to space- and time-limiting procedures which control nonphysical oscillations, makes the implicit time integration difficult, e.g.~because the discrete system is nonlinear also on linear problems. This nonlinearity of the scheme is circumvented as proposed in (Puppo et al., Comm.~Appl.~Math.~\& Comput., 2023) for scalar conservation laws, where a first order implicit predictor is computed to freeze the nonlinear coefficients of the essentially non-oscillatory space reconstruction, and also to assist limiting in time. In addition, we propose a novel conservative flux-centered a-posteriori time-limiting procedure using numerical entropy indicators to detect troubled cells. The numerical tests involve classical and artificially devised stiff problems using the Euler's system of gas-dynamics.
\end{abstract}


\section{Introduction} \label{sec:intro}

Mathematical models for the description of fluids, plasmas, and many other physical phenomena, are typically given in terms of systems of hyperbolic conservation laws. These systems are characterized by a set of multi-dimensional partial differential equations (PDEs) that express the conservation of various physical quantities in terms of their respective fluxes. A prototypical example is provided by the Euler's equations for gas-dynamics describing the conservation of mass, momentum, and energy of a gas.

In this work we focus on one-dimensional systems of $m\geq1$ hyperbolic conservation laws:
\begin{equation} \label{eq:hyp:sys}
	\frac{\partial}{\partial t} \mathbf{u}(x,t) + \frac{\partial}{\partial x} \mathbf{f}(\mathbf{u}(x,t)) = \mathbf{0},
\end{equation}
where, $\mathbf{u}: \mathbb{R} \times \mathbb{R}^+_0 \to \mathbb{R}^m$ is the quantity of interest, and $\mathbf{f}: \mathbb{R}^m \to \mathbb{R}^m$ is the vector of the flux functions. System~\eqref{eq:hyp:sys} is hyperbolic when the eigenvalues $\{\lambda_j(\mathbf{u}(x,t))\}_{j=1}^m$ of the associated Jacobian matrix are real and determine a complete set of eigenvectors. The eigenvalues of~\eqref{eq:hyp:sys} provide the characteristic velocities, which describe the propagation speed of waves in the system. These waves can be either acoustic waves (shocks and rarefactions) or material waves (contact discontinuities). Requiring that the eigenvalues are real implies that the propagation speed of information through the system is finite.

Solving hyperbolic systems of conservation laws is a challenging task, both analytically and numerically, e.g.~due to the occurrence of singularities or the need of devising high order accurate non-oscillatory methods to avoid low-resolution approximations. Another source of numerical difficulty is represented by stiff problems that occur when the system is characterized by speeds spanning different orders of magnitude, namely when $\frac{\max_{j=1,\dots,m} |\lambda_j(\mathbf{u})|}{\min_{j=1,\dots,m} |\lambda_j(\mathbf{u})|} \gg 1$. This happens, for instance, in gas-dynamics when the fluid speed is much less than the speed of the acoustic waves. In many applications the phenomenon of interest travels with a low speed. An example is provided by low-Mach number problems occurring when the equations governing the flow become stiff due to the very low fluid velocity compared to the speed of sound in the fluid. In these situations, the compressibility effects of the fluid can be neglected, and the fluid is almost incompressible. Then, if the interest is on the movement of the fluid, accuracy in the propagation of sound waves becomes irrelevant. For low-Mach problems we refer to~\cite{2010DellacherieLowMach,2011DegondTang,2017AbbateAllSpeed,2017DimarcoLoubereVignal,2018BoscarinoRussoScandurra,2017Tavelli_SemiImplicitAllMach}.

Numerical schemes used to solve hyperbolic problems need to be carefully designed to handle the stiff regime. In fact, it is well-known that explicit schemes are subject to the Courant-Friedrichs-Levy (CFL) stability condition that specifies a constraint on the numerical speed in relation to the maximum speed of information propagating in the system. More precisely, let $\Delta t$ and $h$ be the time-step and the mesh width of a numerical scheme, respectively. We define the \emph{numerical speed} as $s_n = \nicefrac{h}{\Delta t}$ which approximates the speed with which the numerical data propagate in the discretized system. Then, the CFL condition imposes that $s_n$ must be faster than the maximum speed of propagation to ensure that information does not travel too far between adjacent space cells during one time-step. For this reason, the stability request on the time-step of explicit schemes becomes very restrictive for stiff problems due to the presence of fast waves, thus limiting the computational efficiency of the scheme. In contrast, implicit schemes can have superior stability properties. Therefore, they can be less constrained by the CFL condition and can be employed with a time-step focusing on the phenomenon of interest, for instance, on the fluid speed, thus potentially allowing for larger time-step sizes. In addition, since the accuracy of a scheme depends on the difference between the numerical and the actual speed, it turns out that implicit schemes reduce accuracy on the (fast) acoustic waves, still resulting highly accurate on the (slow) material waves. However, implicit methods are computationally more expensive than explicit ones since they require the solution of a system of equations, in general nonlinear, at each time-step. Therefore, the choice between explicit and implicit schemes depends on the specific problem being solved and the desired trade-off between computational efficiency and accuracy.

Here, we deal with an efficient formulation of implicit high order finite volume schemes~\cite{LeVeque:book}. In first order implicit schemes the only source of nonlinearity is due to the nonlinear flux function, namely to the physical structure of the model, which is therefore unavoidable; however they produce large dissipation errors. High order accurate implicit schemes require, as much as their explicit counterparts, nonlinear space-limiting procedures to prevent spurious oscillations (see e.g. \cite{PSV23:Quinpi} for a discussion of the TVD property of implicit schemes). Such space-limiting procedures introduce an additional source of nonlinearity which becomes computationally challenging when using implicit schemes. The novel idea of~\cite{PSV23:Quinpi} was to simplify considerably the implicit high order scheme by using a first order predictor to freeze the non-linearities of the space-limiting procedure. The implicit approach proposed in~\cite{PSV23:Quinpi}, named \emph{Quinpi}, was tailored to the third order implicit approximation of \emph{scalar} conservation laws. More specifically, third order accuracy was achieved by using a third order Diagonally Implicit Runge-Kutta (DIRK) for the time integration and a third order Central Weighted Essentially Non-Oscillatory (CWENO) reconstruction, cf.~\cite{CPSV:cweno}, for the space discretization. The first order implicit scheme was based on a composite backward Euler, evaluated at the abscissae of the DIRK, naturally combined with a piecewise constant (i.e.~linear in the data) reconstruction in space. This predictor was used to freeze the nonlinear weights of the CWENO reconstruction making the resulting third order implicit scheme nonlinear just because of the nonlinearity of the flux function.

As noted in~\cite{2020Arbogast,PSV23:Quinpi}, the appearance of spurious oscillations can be observed in implicit integration, especially for large Courant numbers, despite space-limiting being performed. A \emph{time-limiting} procedure is required in order to control the oscillations arising in the implicit time integration, i.e.~when the large time-step size allows for propagation of waves crossing several adjacent space cells during a single time-step. In this context, ``time-limiting'' refers to the practice of downgrading the order of the solution on irregular cells, thereby mitigating spurious oscillations associated with the use of large time-steps. The problem of the time-limiting has been discussed in several papers, which typically deal with limiting in space
and time simultaneously. For instance, we mention~\cite{2003DurasaisamyBaeder,2007DurasaisamyBaeder} for second-order schemes and~\cite{2020Arbogast} for a fully nonlinear third order implicit scheme. In~\cite{PSV23:Quinpi}, instead, the time-limiting was obtained as a cell-centered a-posteriori nonlinear blending of the first order and third order solutions. Troubled cells were detected using a combination of space and time regularity indicators, in a WENO-like fashion. However, being cell-centered, the technique in~\cite{PSV23:Quinpi} has the drawback of being non-conservative, and, thus, it requires a conservative correction. Conservative flux-centered time-limiting procedures, inspired by the Multi-dimensional Optimal Order Detection (MOOD) method introduced in~\cite{CDL11:MOOD,CDL12:MOOD}, were proposed in~\cite{VTSP23:Quinpi:Book,EBS22:Implicit:Networks} for the control of nonphysical oscillations of high order implicit numerical solutions. MOOD was originally developed in order to reduce the order of the space reconstruction on problematic cells with the help of several problem-dependent detectors to check whether extrema of the numerical solution are smooth, physical or spurious oscillations. MOOD was also extended to other contexts, as for instance in~\cite{SL:18:AMRMOOD,LDD:14,ZDLS:14}.

The contribution of the present work is twofold:
\begin{enumerate}
	\item We extend the Quinpi framework proposed in~\cite{PSV23:Quinpi} to the case of the numerical approximation of general stiff hyperbolic systems when the interest is on the movement of the fluid. This implicit approach is not tailored to the solution of a specific stiff problem, as it usually happens when devising numerical schemes for low-Mach problems. Moreover, the space-limiting exploits the novel CWENO reconstruction~\cite{STP23:cweno:boundary} which does not make use of ghost cells for boundary reconstructions;
	\item We introduce a flux-centered conservative time-limiting procedure inspired by the MOOD technique, namely we replace the high order numerical fluxes at the cell interfaces of troubled cells with low order numerical fluxes. A crucial point is the choice of the troubled cells indicator; here we investigate the use of the numerical entropy production error~\cite{PS11:numerical:entropy,Puppo04:numerical:entropy} as indicator that signals non-smooth solutions instead of the typical MOOD detectors of oscillatory cells. 
\end{enumerate}

Finally, we mention the following approaches to the implicit integration of hyperbolic systems developed in the literature. A fifth order implicit WENO scheme was proposed in~\cite{2001Gottlieb}, where a predictor-corrector technique was also used. However, the predictor was based on an explicit first order scheme and therefore it allows to deal with small Courant numbers only. A fully nonlinear implicit scheme, based on a third order RADAU time integrator and a third order WENO reconstruction, can be found in~\cite{2020Arbogast}. Fully implicit, semi-implicit, implicit-explicit, local time-stepping and active flux treatments of stiff hyperbolic equations were also investigated, e.g., in~\cite{CNPT:09,CNPT:10,CNPT:10a,CPPT:06,FZ:22,FKRZ:22,ABIR:19,BB:23,BoscarinoQiuRussoXiong}.

The paper is organized as follows. In Section~\ref{sec:scheme} we introduce the fully implicit third order scheme of system~\eqref{eq:hyp:sys} which uses the CWENO reconstruction of~\cite{STP23:cweno:boundary} for the space approximation and a DIRK method for the time integration. Then, in Section~\ref{sec:quinpi}, we devise the Quinpi framework, which allows to overcome the nonlinearity of CWENO by freezing the nonlinear weights with a first order composite backward Euler approximation of the solution. The limiting in time is discussed, which employs a flux-centered MOOD technique combined with the numerical entropy production error as detector of troubling cells. Numerical experiments are provided in Section~\ref{sec:numerics}. Precisely, Section~\ref{sec:numerics:euler} and Section~\ref{sec:numerics:lowMach} are focused on the Euler's equations for gas-dynamics, where we investigate the experimental order of convergence of the Quinpi scheme and its performance on stiff and low-Mach problems. Section~\ref{sec:numerics:scalar} is devoted to the comparison of the novel time-limiting procedure presented in this work, the cell-centered one of~\cite{PSV23:Quinpi} and the flux-centered of~\cite{VTSP23:Quinpi:Book}. Finally, we discuss results and perspectives in Section~\ref{sec:conclusion}.

\section{Third order space-time fully implicit discretization} \label{sec:scheme}

We consider a finite volume approximation of~\eqref{eq:hyp:sys} through the method of lines (MOL) on the compact computational domain $\Omega=[a,b]\subset\mathbb{R}$. To this end, we discretize $\Omega$ with $N$ uniform cells $\Omega_j=[x_j-\nicefrac{h}{2},x_j+\nicefrac{h}{2}]$ of amplitude $h>0$, such that $\cup_{j=1}^N \Omega_j=\Omega$ and $x_j = a + (j-\nicefrac12)h$ are the cell centers. For the sake of simplicity, we will describe the scheme for a uniform grid, but it is easy to generalize the scheme to a non uniform mesh. 
Defining the cell averages of the exact solution on a given space cell $\Omega_j$ as
$$
\ca{\mathbf{u}}_j(t)= \frac1h \int_{\Omega_j} \mathbf{u}(x,t) \dx, \quad t \geq 0
$$
the MOL provides the following semi-discrete form of system~\eqref{eq:hyp:sys}:
\begin{equation} \label{eq:MOL}
	\TotDer{\ca{\mathbf{u}}_j(t)}{t} = -\frac{1}{h}\left[ \mathbf{f}\left(\mathbf{u}_{j+\frac12}(t)\right) - \mathbf{f}\left(\mathbf{u}_{j-\frac12}(t)\right)\right], \quad j=1,\dots,N, \ t\geq 0,
\end{equation}
where $\mathbf{u}_{j\pm\nicefrac12}(t) = \mathbf{u}\left(x_j\pm\nicefrac{h}{2},t\right)$. System~\eqref{eq:MOL} describes the conservation of the cell averages as the difference of the right and left fluxes at the cell boundaries. Up to now, no numerical approximation of the exact solution has been introduced. In fact, \eqref{eq:MOL} is still exact.

In order to transform the MOL in a numerical scheme, first one has to transform the exact system~\eqref{eq:MOL} in a closed system for the cell averages and, afterwards, one can introduce a time discretization of the resulting coupled ODE system which evolves in time the cell averages. Consequently, there is the problem of the knowledge of point values of the solution at the cell interfaces for the evaluation of the flux function. The extrapolation of these values from the cell averages is the so-called \emph{reconstruction problem}. WENO and CWENO schemes, see e.g.~\cite{JiangShu:96,Shu:97,LPR:00:SIAMJSciComp,CPSV:cweno}, and their developments~\cite{Balsara:AOWENO,CCD:11,ABC16:improvedWENOZ,SempliceVisconti:2020}, are examples of numerical procedures computing point values as function of the cell averages. The advantage of WENO and CWENO schemes is that they achieve high order approximations of the reconstructions, but they pay the price of enlarging the stencil, compared to low order schemes, and of being highly nonlinear, which is a computational bottleneck in implicit time integration.

In the following, before dealing with the time integration of~\eqref{eq:MOL}, we recall the space reconstruction based on a recently developed CWENO scheme~\cite{STP23:cweno:boundary}.

\subsection{Space reconstruction: third order CWENOZ without ghost cells} \label{sec:scheme:cweno}

In reconstruction procedures, the goal is to provide a space limited approximation of the exact solution $\mathbf{u}(\cdot,t)$, at a given time $t\geq 0$, using the knowledge of its cell averages. Since the reconstruction is typically applied component-wise, we will describe the procedure on a component $u$ of the vector solution $\mathbf{u}$ in one-dimension.

A CWENO type reconstruction defines an approximation of $u(\cdot,t)$ as
$$
u(x,t) \approx \sum_{j=1}^N R_j(x;t)\chi_{\Omega_j}(x), \quad t\geq 0,
$$
where $\chi_{\Omega_j}$ is the characteristic function of the cell $\Omega_j$, and $R_j(x;t)$ is the space reconstruction polynomial for $x\in\Omega_j$, which depends on time through the time-dependent cell averages. Indeed, here the time variable $t$ is treated as a parameter meaning that the space reconstruction polynomial is computed at a fixed time. If the desired reconstruction point, say $\hat{x}$, lies within the cell $\Omega_j$, then the evaluation of the polynomial $R_j(\hat{x};t)$ provides the needed point values of $u(\hat{x},t)$. The CWENO type procedure differs from the classical WENO scheme by the fact that each polynomial $R_j(x;t)$ is globally defined in its reference cell $\Omega_j$ and, therefore, it can be pre-computed and later evaluated at the needed locations.

We focus on third order space reconstructions. Then, the CWENO scheme with Z-type nonlinear weights~\cite{CCD:11,CSV19:cwenoz} defines the polynomial $R_j(x;t)$ for $x\in\Omega_j$ as follows.

\begin{definition}[Third order CWENOZ reconstruction, see~\cite{CSV19:cwenoz}] \label{def:CWENOZ}
	Let $\Popt\in\Poly{2}$ be the \emph{optimal} polynomial of degree $2$, which interpolates all the data in the three-cell stencil $\Sopt = \{ \Omega_{j-1}, \Omega_j, \Omega_{j+1} \}$. Further, let ${P}_L, {P}_R\in\Poly{1}$ be polynomials of degree $1$ such that $P_L$ interpolates the cell averages of the left-biased sub-stencil $\mathcal{S}_L = \{ \Omega_{j-1}, \Omega_j \}$, and $P_R$ interpolates the cell averages of the right-biased sub-stencil $\mathcal{S}_R = \{ \Omega_j, \Omega_{j+1} \}$.
	Let also $\{d_0,d_L,d_R\}$ be a set of strictly positive real coefficients such that $\sum_{k=0,L,R} d_k=1$.
	
	The CWENOZ procedure computes the reconstruction polynomial on $\Omega_j$ as
	\begin{equation}
		R_j^{\mathsf{CWZ}}(x;t) 
		= \frac{\omega_0}{d_0} \left( \Popt(x;t) - \sum_{k=L,R} d_k P_k(x;t) \right) + \sum_{k=L,R} \omega_k P_k(x;t) \in\Poly{2}, \label{eq:precCWZ}
	\end{equation}
	where $\omega_0$, $\omega_L$ and $\omega_R$ are the (nonlinear) coefficients defined as
	\begin{equation} \label{eq:omegaZ}
		\alpha_k = {d_k} \left( 1 + \left( \frac {\tau} {I_k+\epsilon} \right)^{p} \right),
		\qquad
		\omega_k = \frac{\alpha_k}{\sum_{i=L,0,R}\alpha_i}, \quad k=0,L,R.
	\end{equation}
	In~\eqref{eq:omegaZ}, $I_0$, $I_L$ and $I_R$ are the regularity indicators of the associated polynomials $P_L$, $\Popt$ and $P_R$, respectively, computed as the Jiang-Shu indicators from~\cite{JiangShu:96}:
	\begin{equation} \label{eq:ind}
		I[P_k] := 
		\sum_{i=1}^{\deg(P_k)}  h^{2i-1} \int_{\Omega_j} \left(\frac{\d^i}{\dx^i} P_k(x;t)\right)^2 \dx, \quad k=0,L,R.
	\end{equation}		
	Finally, $\epsilon=h^{q}$, $q\geq 1$, $p \ge 1$ and $\tau$ is the following global smoothness indicator 
	\begin{equation}  \label{eq:tau}
		\tau := \left| 2I_0 - I_L - I_R \right|.
	\end{equation}
\end{definition}

The CWENOZ reconstruction polynomial switches between the high accurate polynomial $\Popt$, when the cell averages in the stencil $\Sopt$ are a sampling of a smooth enough function, and a nonlinear blending of $\Popt$ and of the lower degree polynomials $P_L, P_R$ when a discontinuity is present in the stencil $\Sopt$. 
%The switch is automatically performed thanks to the definition of the nonlinear weights $\omega_0$, $\omega_L$ and $\omega_R$, computed with the help of the regularity indicators~\eqref{eq:ind}. In practice, these are of order $o(1)$ when the stencil $\Sopt$ is smooth, so that $\omega_i \approx d_i$, $i=0,L,R$, and thus $R_j^{\mathsf{CWZ}} \approx \Popt$. In this case, the reconstruction is able to achieve the maximal desired order of accuracy. Instead, the regularity indicator of a polynomial interpolating discontinuous data is of order $\mathcal{O}(1)$, so that the corresponding nonlinear weight deviates from its optimal value, avoiding the appearance of spurious oscillations in the reconstruction polynomial.

The use of the Z-type weights~\eqref{eq:omegaZ} allows to have better accuracy on smooth data compared to classical weights~\cite{CPSV:cweno,Shu:97,LPR:00:SIAMJSciComp}, especially on coarse grids, without sacrificing the non-oscillatory properties. This is obtained by using the optimal choice~\eqref{eq:tau} of the global smoothness indicator that makes $\tau$ much smaller than the regularity indicators when the data in $\Sopt$ are smooth enough. %Typically, $\tau$ is restricted to be a linear combination of the other smoothness indicators for efficiency. For results on the optimal choices of $\tau$ we refer to~\cite{CSV19:cwenoz}.
In~\cite{CSV19:cwenoz}, it is proven that the accuracy of the CWENOZ reconstruction on smooth flows is the optimal one, provided that $\deg(\Popt) \leq 2\deg(P_k)$, for $k = L,R$. %This condition is satisfied by the third-order reconstruction in Definition~\ref{def:CWENOZ}.

Close to boundaries, the central stencil $\Sopt$ may not be defined, because it would not be fully contained in $\Omega$. In this case, one can consider for $\Popt$ a 3-cell stencil entirely biased towards the domain interior, but the need of controlling spurious oscillations requires the inclusion of a polynomial $\tilde{P}\in\Poly{0}$, defined on the cells which contain the endpoints of the physical domain. Indeed, assume that the 3-cell stencil is given by $\{\Omega_1,\Omega_2,\Omega_3\}$, then $\tilde{P}$ allows to select the smooth part of the stencil when a discontinuity is present either in $\Omega_2$ or $\Omega_3$. Optimal accuracy can still be achieved provided that the corresponding linear weight is infinitesimal of order $\mathcal{O}(h^r)$, for some $r > 0$. This approach was introduced for CWENOZ type reconstructions in~\cite{SempliceVisconti:2020}, where a thorough study of sufficient conditions on $r$, and on the other parameters of the scheme, to achieve optimal accuracy has been performed, and exploited for reconstructions free of ghost cells in \cite{STP23:cweno:boundary}. In the following definition, we recall a particular third order CWENOZ Adaptive Order reconstruction which will be used to define the boundary reconstruction. Again, we consider the one-dimensional case. For two-dimensional reconstructions we refer to~\cite{STP23:cweno:boundary}.
%can be recast in the framework of~\cite{SempliceVisconti:2020}.

\begin{definition}[Third order CWENOZ reconstruction without ghost cells, see~\cite{STP23:cweno:boundary}] \label{def:CWENOZb}
	Let $R_j(x;t)$ be the reconstruction polynomial related to the cell $\Omega_j$. Then,
	\begin{equation} \label{eq:precCWZb}
		R_j(x;t) =
		\begin{cases}
			R_j^{\mathsf{CWZ}}(x;t), & j=2,\dots,N-1, \\[1ex]
			R_j^{\mathsf{AO}}(x;t), & j=1 \ \mbox{ with } \ \Sopt=\{\Omega_1,\Omega_2,\Omega_3\}, \\[1ex]
			R_j^{\mathsf{AO}}(x;t), & j=N \ \mbox{ with } \ \Sopt=\{\Omega_{N-2},\Omega_{N-1},\Omega_N\},
		\end{cases}
	\end{equation}
	where $\tau_j = \tau$ in~\eqref{eq:tau}, for $j=2,\dots,N-1$, and $\tau_1=\tau_2$, $\tau_N=\tau_{N-1}$.
\end{definition}

The reconstruction in the first and last computational cells appearing in \eqref{eq:precCWZb} are defined exploting the results  in \cite{SempliceVisconti:2020}, which we report here for completeness in the third order case.

\begin{definition}[Third order CWENO-AO reconstruction, see~\cite{SempliceVisconti:2020}] \label{def:CWENOZAO}
	Let $\Popt\in\Poly{2}$ be the \emph{optimal} polynomial of degree $2$ which interpolates all the given data in the three-cell stencil $\Sopt$ such that $\Omega_j \in \Sopt$. Let $P\in\Poly{1}$ be the polynomial of degree $1$ such that $P$ interpolates the cell averages of a two-cell sub-stencil $\mathcal{S}$ such that $\Omega_j\in\mathcal{S}\subset\Sopt$. Further, let $\tilde{P}\in\Poly{0}$ be the constant polynomial associated to the sub-stencil $\tilde{\mathcal{S}} = \{ \Omega_j \}$, namely $\tilde{P}(x;t)=\bar{u}_j$ for $x\in\Omega_j$.
	Let also $\{d_0,d,\tilde{d}\}$ be a set of strictly positive real coefficients such that $d_0+d+\tilde{d}=1$ and $\tilde{d}=h^r$ for some $r>0$.
	
	The CWENOZ-AO reconstruction polynomial on $\Omega_j$ is
	\begin{equation}
		R_j^{\mathsf{AO}}(x;t) 
		= \frac{\omega_0}{d_0} \left( \Popt(x;t) - d P(x;t) - \tilde{d} \tilde{P}(t) \right) + \omega P(x;t) + \tilde{\omega} \tilde{P}(t) \in\Poly{2}, \label{eq:precAO}
	\end{equation}
	where $\omega_0$, $\omega$ and $\tilde{\omega}$ are the (nonlinear) Z-type coefficients, see~\eqref{eq:omegaZ}, associated to the regularity indicators $I_0$, $I$ and $\tilde{I}$, respectively for the polynomials $\Popt$, $P$ and $\tilde{P}$, see~\eqref{eq:ind}.
	
	Since $\tilde{I}=0$, the global smoothness indicator $\tau$ 
	in this case can be chosen as
	\begin{equation} \label{eq:tauAO}
		\tau := \left| I -  I_0 \right|,
	\end{equation}
\end{definition}

The analysis in~\cite{SempliceVisconti:2020} shows that the reconstruction of Definition~\ref{def:CWENOZb} achieves third order of accuracy for $r=1,2$, provided that the exponent $p$ in \eqref{eq:omegaZ} is $p \geq 1$ and that $\epsilon=h^q$ for $q=1,2,3$.

%\begin{notation} \label{notation:poly}
%In the following, for convenience of notation, we denote $P_L$ the constant polynomial on its stencil denoted by $\mathcal{S}_L=\{\Omega_1\}$ and $P_R$ the linear polynomial interpolating the data in its stencil denoted by $\mathcal{S}_R=\{\Omega_1,\Omega_2\}$ when the reconstruction on $\Omega_1$ is performed. Symmetrically, for the last cell.
%\end{notation}

%In the following, for convenience, we illustrate the boundary reconstruction on the left boundary for the cell $\Omega_1$.

The reconstruction polynomial $R_j(x;t)$ given in~\eqref{eq:precCWZb} provides the approximation for each component $u_j$ of $\mathbf{u}_j$ for $x\in\Omega_j$. Thus, one can estimate the values $u\left(x_j\pm\nicefrac{h}{2},t\right)$, $t\geq 0$, with
\begin{subequations} \label{eq:bed}
	\begin{gather}
		u^-_{j+\frac12}(t)=R_{j}\left(x_j+\frac{h}{2};t\right) \ \text{ and } \ u^+_{j+\frac12}(t)=R_{j+1}\left(x_j+\frac{h}{2};t\right), \quad j=1,\dots,N-1 \label{eq:bed:interior} \\
		u^-_{\frac12}(t)=u^-_{\mathsf{out}} \ \text{ and } \ 	u^+_{\frac12}(t)=R_1\left(a;t\right) \label{eq:bed:first} \\
		u^-_{N+\frac12}(t)=R_N\left(b;t\right) \ \text{ and } \ 	u^+_{N+\frac12}(t)=u^+_{\mathsf{out}} \label{eq:bed:last}	
	\end{gather}
\end{subequations}
which are named \emph{boundary extrapolated data} (BED). The outer values $u^{\mp}_{\mathsf{out}}$ are determined by the boundary conditions. For instance, for periodic boundary conditions the outer values $u^{\mp}_{\mathsf{out}}$ are set to the inner reconstructions at the last and first interfaces, respectively. Instead, for free-flow boundary conditions, the outer values $u^{\mp}_{\mathsf{out}}$ are set to the inner reconstructions at the first and last interfaces, respectively, namely imposing zero jumps at the boundaries of the physical domain.

Notice that at each interface the two BED in~\eqref{eq:bed} are different, although computed at the same interface $x_{j+\nicefrac{h}{2}}$. Therefore, in order to approximate the flux function at the interfaces, one introduces a consistent and monotone numerical flux function
\begin{subequations} \label{eq:numfluxfnc}
	\begin{equation}
		(\mathbf{v},\mathbf{w}) \in \R^m\times\R^m \mapsto \NumFlux(\mathbf{v},\mathbf{w}) \in \R^m,
	\end{equation}
	such that
	\begin{equation}
		\mathbf{f}\left(\mathbf{u}(x_{j+\frac{1}{2}},t)\right) \approx %\NumFlux_{j+\nicefrac12}(t) = 
		\NumFlux\left(\mathbf{u}^-_{j+\frac12}(t),\mathbf{u}^+_{j+\frac12}(t)\right) \in \R^m.
	\end{equation}
\end{subequations}
The function $\NumFlux$ may be any approximate or exact Riemann solver, which is applied component-wise on vector-valued inputs. Finally, the exact system of ODEs~\eqref{eq:MOL} is reduced to a finite system of ODEs for the evolution of the cell averages. The right hand side is completely defined by the space reconstruction along with the numerical flux, and one obtains
\begin{equation} \label{eq:spaceApprox}
	\TotDer{\ca{\mathbf{U}}_j(t)}{t} = - \frac{1}{h}\left[ \NumFlux_{j+\frac12}(t)-\NumFlux_{j-\frac12}(t)\right],
\end{equation}
which provides the approximation $\ca{\mathbf{U}}_j(t)$ of the cell averages $\ca{\mathbf{u}}_j(t)$ of the solution $\mathbf{u}(x,t)$, $x\in\Omega_j$, and where
\begin{equation} \label{eq:numflux}
	\NumFlux_{j+\frac12}(t) = \NumFlux\left(\mathbf{U}^-_{j+\frac12}(t),\mathbf{U}^+_{j+\frac12}(t)\right) \in \R^m,
\end{equation}
with $\mathbf{U}^{\mp}_{j+\nicefrac12}(t)$ BED of the data $\ca{\mathbf{U}}(t) = \left[\ca{\mathbf{U}}_1(t),\dots,\ca{\mathbf{U}}_N(t)\right]^T$ according to~\eqref{eq:bed}.

%In the following, we use the reconstruction of Definition~\ref{def:CWENOZb} in the implicit time treatment of system~\eqref{eq:spaceApprox}.

\subsection{Time integration: third order Diagonally Implicit Runge-Kutta} \label{sec:scheme:dirk}

In order to employ a timestep $\Delta t$ which is not constrained by the CFL stability condition,
we solve numerically equation~\eqref{eq:spaceApprox} with a Diagonally Implicit Runge-Kutta (DIRK) scheme with $s$ stages and general Butcher tableau
\begin{equation}
	\label{eq:tableau:dirk}
	\begin{array}{c|cccc}
		c_1 & a_{11} & 0 & \dots & 0 \\[1.5ex]
		c_2 & a_{21} & a_{22} & \dots & 0 \\[1.5ex]
		\vdots & \vdots & \vdots & \ddots & \\[1.5ex]
		c_s & a_{s1} & a_{s2} & \dots & a_{ss} \\[1ex]
		\hline
		&&&&\\[-1.8ex]
		& b_1 & b_2 & \dots & b_s
	\end{array}
\end{equation}
A typical assumption is that $c_k = \sum_{i=1}^s a_{ki}$, $k=1,\dots,s$, and one has $\sum_{k=1}^s b_k=1$ for consistency. 
Further, choosing a scheme in which $a_{kk}$ is independent of $k$, the construction of the Jacobian in the nonlinear solvers is simplified.
Clearly, one uses a Butcher tableau with order matching the order of the space reconstruction, in this case a third order accurate scheme. Therefore, DIRK schemes with number of stages $s\geq 2$ must be considered, e.g.~see~\cite{Alexander1977}.

The space-time discretization leads to the fully-discrete scheme
\begin{subequations} \label{eq:implicit}
	\begin{align}
		\ca{\mathbf{U}}_j^{(k)} &= \ca{\mathbf{U}}_j^{n} - \frac{\Delta t}{h} \sum_{i=1}^k a_{ki} \left[ \NumFlux_{j+\frac12}^{(i)} - \NumFlux_{j-\frac12}^{(i)} \right], \quad k=1,\dots,s, \label{eq:implicit:stage} \\
		\ca{\mathbf{U}}_j^{n+1} &= \ca{\mathbf{U}}_j^{n} - \frac{\Delta t}{h} \sum_{k=1}^s b_k \left[ \NumFlux_{j+\frac12}^{(k)} - \NumFlux_{j-\frac12}^{(k)} \right], \quad n\geq 0, \label{eq:implicit:update} \\
		\NumFlux_{j+\frac12}^{(k)} &= \NumFlux\left(\mathbf{U}_{j+\frac12}^{-,(k)},\mathbf{U}_{j+\frac12}^{+,(k)}\right) \in \R^m \label{eq:implicit:flux}
	\end{align}
\end{subequations}
for each $j=1,\dots,N$, where $\Delta t$ is the time-step, $\ca{\mathbf{U}}_j^{n} \approx \ca{\mathbf{u}}_j(n\Delta t)$, and $\mathbf{U}_{j+\nicefrac12}^{\mp,(k)}$ are the BED of the stage values $\ca{\bigvec{U}}^{(k)} = \left[ \ca{\mathbf{U}}_1^{(k)}, \dots, \ca{\mathbf{U}}_N^{(k)} \right]^T$, 
approximations of the solution at times $t^{(k)} = (n + c_k) \Delta t$, according to~\eqref{eq:bed}. Here and in the following, the notation suggests the use of a uniform time-step $\Delta t$. Nevertheless, the scheme can be formulated for a non-uniform time-step, but we prefer not to burden the notation.

The advantage of DIRK schemes is that the implicit computation of a given stage value~\eqref{eq:implicit:stage} can be performed sequentially from $k=1$ to $k=s$. Therefore, at each time-step one has to solve $s$ systems of nonlinear equations of size $mN$ by $mN$:
\begin{equation}\label{eq:DIRK:stage}
	\mathbf{G}\left(\ca{\bigvec{U}}^{(k)}\right) := \ca{\bigvec{U}}^{(k)} 
	+ \frac{a_{kk}\Delta t}{h} \Delta \bigvec{F}^{(k)} 
	- \ca{\bigvec{U}}^{n} 
	+ \frac{\Delta t}{h} \sum_{i=1}^{k-1} a_{ki} \Delta \bigvec{F}^{(i)} = \mathbf{0},
\end{equation}
where $\Delta\bigvec{F}^{(k)} \in \R^{mN}$ whose $j$-th block is 
$ \Delta\bigvec{F}^{(k)}_j = \left( \NumFlux_{j+\nicefrac12}^{(k)}-\NumFlux_{j-\nicefrac12}^{(k)} \right) \in \R^m$, for $j=1,\dots,N$.
Above we have highlighted the term $\Delta \bigvec{F}^{(k)}$, which makes the system for the $k$-th stage value nonlinear. In fact, the other flux differences $\Delta \bigvec{F}^{(i)}$, $i=1,\dots,k-1$, are already available, thanks to the structure of DIRK schemes.

As already noticed in~\cite{PSV23:Quinpi}, $\mathbf{G}$ has two sources of nonlinearity. One is unavoidable because it is due to the physics when the phenomenon under study is described by the nonlinear flux function $\mathbf{f}$ in~\eqref{eq:hyp:sys}. The second, instead, is introduced by the high order space reconstruction procedure, needed for the computation of the BED, which is highly nonlinear because of the nonlinear weights~\eqref{eq:omegaZ} and of the regularity indicators~\eqref{eq:ind}. Therefore, even for linear PDEs, a standard implicit scheme requires a nonlinear solver to find the solution of $\mathbf{G}(\ca{\bigvec{U}}^{(k)})=\mathbf{0}$, in equation \eqref{eq:DIRK:stage} for $k=1,\dots,s$. Typically, one uses Newton-Raphson's method, which requires assembling the Jacobian of the nonlinear system $\mathbf{G}$, resulting in a high computational cost. In fact, the Jacobian required for the Newton iterations has bandwidth one point larger than the stencil size in each direction, so for a third order scheme it is a block-pentadiagonal matrix with $m\times m$ blocks, with entries depending on the nonlinear weights and on the regularity indicators, which have very complicated expressions.

The Quinpi approach, introduced in~\cite{PSV23:Quinpi} for scalar conservation laws, provides a way to circumvent the nonlinearity determined by the high order reconstruction procedure. In the following section, we extend the Quinpi idea to the hyperbolic system~\eqref{eq:hyp:sys}.

\section{Third order Quinpi scheme for one-dimensional hyperbolic systems}
\label{sec:quinpi}

The name Quinpi stands for implicit CWENO and it is based on a predictor-corrector approach to avoid the nonlinearity of the high order scheme, introduced by the space reconstruction, and keep the nonlinearity of the flux function $\mathbf{f}$ only.

The Quinpi idea relies on the following considerations on the structure of most essentially-non-oscillatory reconstructions. Since reconstructions are typically applied component-wise, we use the notation for scalar conservation laws in this section.

Let $\mathcal{S}_p = \{ \Omega_{j-p}, \dots, \Omega_{j+p} \}$ be a stencil of $2p+1$ computational cells, and $P\in\Poly{2p}$ be the optimal interpolating polynomial of the cell averages in $\mathcal{S}_p$. Then, the dependence of $P$ on the data is linear, i.e.
\begin{equation*}
	P(x;t) = \sum_{\alpha=-p}^{p} \mu_{j,\alpha}(x) \ca{u}_{j+\alpha}(t).
\end{equation*}
Therefore, unrolling the linearity with respect to the data of all the interpolating polynomials involved in the third order reconstruction procedure of Definition~\ref{def:CWENOZb}, one can write the reconstruction polynomial as
\begin{equation*}
	\begin{aligned}
		R_j(x;t) &= \sum_{i=0,L,R} \omega_i\left(\{\ca{u}_k(t)\}_{k\in\Sopt}\right)	P^i_j(x) \\
		&= \sum_{i=0,L,R} \omega_i\left(\{\ca{u}_k(t)\}_{k\in\Sopt}\right)	\sum_{\alpha=-1}^{1} \mu^{i}_{j,\alpha}(x) \ca{u}_{j+\alpha+\delta_{j1}-\delta_{jN}}(t) \\
		&= \sum_{\alpha=-1}^{1} W_{j,\alpha} \left(x;\{\ca{u}_k(t)\}_{k\in\Sopt}\right)  \ca{u}_{j+\alpha+\delta_{j1}-\delta_{jN}}(t),
	\end{aligned}
\end{equation*}
where we have collected the nonlinear weigths in the quantities
\begin{equation*}
	W_{j,\alpha} \left(x;\{\ca{u}_k(t)\}_{k\in\Sopt}\right)
	= \sum_{i=0,L,R} \omega_i\left(\{\ca{u}_k(t)\}_{k\in\Sopt}\right)	\mu^{i}_{j,\alpha}(x),
\end{equation*}
and $\delta_{ij}$ denotes the Kronecker delta.
For the first and the last cell, the range of summation over $\alpha$ is adjusted according to the stencil of the reconstruction.
%where we have used the notation convention~\ref{notation:poly}.
Then, it is possible to write the inner BED in~\eqref{eq:bed} as
%\begin{subequations} \label{eq:bed:split}
\begin{equation} \label{eq:bed:split}%\label{eq:bed:split:inner}
	\begin{aligned}
		u_{j+\frac12}^{-}(t) &= R_j\left(x_{j+\frac12};t\right) \\
		&= \sum_{\alpha=-1}^1 W_{j,\alpha}\left(x_{j+\frac12};\{\ca{u}_k(t)\}_{k\in\Sopt}\right) \overline{u}_{j+\alpha+\delta_{j1}-\delta_{jN}}(t), \\
		u_{j+\frac12}^{+}(t) &= R_{j+1}\left(x_{j+\frac12};t\right) \\
		&= \sum_{\alpha=-1}^1 W_{j+1,\alpha}\left(x_{j+\frac12};\{\ca{u}_k(t)\}_{k\in\Sopt}\right) \overline{u}_{j+1+\alpha+\delta_{j1}-\delta_{jN}}(t),
	\end{aligned}
\end{equation}
% for $j=1,\dots,N-1$, and
% 	\begin{equation} \label{eq:bed:split:bound}
	% 		\begin{aligned}
		% 		u_{\frac12}^{+}(t) &= R_1\left(a;t\right) = \sum_{\alpha=-1}^1 W_{1,\alpha}\left(a;\{\ca{u}_k(t)\}_{k\in\Sopt}\right) \overline{u}_{2+\alpha}(t), \\
		% 		%
		% 		u_{N+\frac12}^{-}(t) &= R_{N}\left(b;t\right) = \sum_{\alpha=-1}^1 W_{N,\alpha}\left(b;\{\ca{u}_k(t)\}_{k\in\Sopt}\right) \overline{u}_{N+\alpha-1}(t).
		% 		\end{aligned}
	% 	\end{equation}
%\end{subequations}
We have highlighted the dependence of $W_{j,\alpha}$ on the data $\{\ca{u}_k\}_{k\in\Sopt}$, since this is highly nonlinear because it contains the nonlinear weights.

\begin{remark}
	For the third order reconstruction considered in this work, the quantities $W_{j,\alpha} \left(x;\{\ca{u}_k(t)\}_{k\in\Sopt}\right)$ have the following expressions. For the left BED at $x=x_{j+\nicefrac12}$, with $j=2,\dots,N-1$:
	\begin{gather*}
		W_{j,-1} = \left( \frac{-1+3 d_L}{6} \right) \frac{\omega_0}{d_0} - \frac12 \omega_L, \quad W_{j,0} = \left( \frac{5-9d_L-3d_R}{6} \right) \frac{\omega_0}{d_0} + \frac32 \omega_L + \frac12 \omega_R \\
		W_{j,1} = \left( \frac{2-3 d_R}{6} \right) \frac{\omega_0}{d_0} + \frac12 \omega_R.
	\end{gather*}
	For the right BED at $x=x_{j+\nicefrac12}$, with $j=1,\dots,N-2$:
	\begin{gather*}
		W_{j+1,-1} = \left( \frac{2-3 d_L}{6} \right) \frac{\omega_0}{d_0} + \frac12 \omega_L, \quad W_{j+1,0} = \left( \frac{5-3d_L-9d_R}{6} \right) \frac{\omega_0}{d_0} + \frac12 \omega_L + \frac32 \omega_R \\
		W_{j+1,1} = \left( \frac{-1+3 d_R}{6} \right) \frac{\omega_0}{d_0} - \frac12 \omega_R.
	\end{gather*}
\end{remark}

The main idea of Quinpi is to exploit a predictor
$\{\ca{u}^{\star}_j(t)\}_{j=1}^N$
of the solution $\{\ca{u}_j(t)\}_{j=1}^N$ at time $t$ and use it to pre-compute and freeze the nonlinear weights in~\eqref{eq:bed:split}, so that the BED can be approximated as
\begin{equation} \label{eq:bed:linearized}
	\begin{aligned}
		u_{j+\frac12}^{-}(t) &\approx \hat{u}_{j+\frac12}^{-}(t) := \sum_{\alpha=-1}^1 W_{j,\alpha}\left(x_{j+\frac12};\{\ca{u}^{\star}_k(t)\}_{k\in\Sopt}\right) \overline{u}_{j+\alpha+\delta_{j1}-\delta_{jN}}(t), \\
		u_{j+\frac12}^{+}(t) &\approx \hat{u}_{j+\frac12}^{+}(t) := \sum_{\alpha=-1}^1 W_{j+1,\alpha}\left(x_{j+\frac12};\{\ca{u}^{\star}_k(t)\}_{k\in\Sopt}\right) \overline{u}_{j+1+\alpha+\delta_{j1}-\delta_{jN}}(t). \\
		% u_{\frac12}^{+}(t) &\approx \hat{u}_{\frac12}^{+}(t) := \sum_{\alpha=-1}^1 W_{1,\alpha}\left(a;\{\ca{u}^{\star}_k(t)\}_{k\in\Sopt}\right) \overline{u}_{2+\alpha}(t), \\
		% u_{N+\frac12}^{-}(t) &\approx \hat{u}_{N+\frac12}^{-}(t) := \sum_{\alpha=-1}^1 W_{N,\alpha}\left(b;\{\ca{u}^{\star}_k(t)\}_{k\in\Sopt}\right) \overline{u}_{N+\alpha-1}(t).
	\end{aligned}
\end{equation}
In this way, the complete scheme would be linear with respect to the space reconstruction, and nonlinear only through the flux function.

The proposal of a predictor to linearize an implicit high order scheme was first discussed in~\cite{Gottlieb:iWENO:2006} for WENO reconstructions. There, the solution of a first order explicit scheme is used as predictor in order to compute the nonlinear WENO weights.
Instead, in~\cite{PSV23:Quinpi,2020Arbogast} an implicit first order scheme was chosen as predictor. Although an implicit predictor is more expensive than the explicit approach proposed in~\cite{Gottlieb:iWENO:2006}, it is stable and non-oscillatory even for higher Courant numbers, thus allowing for reliable prediction of the nonlinear weights.

In this work, we follow the approach in~\cite{PSV23:Quinpi} and describe it for systems of conservation laws below. With respect to the above description, in the case of a system of $m$ conservation laws, one has a reconstruction polynomial $\mathbf{R}_j(x;t)\in(\mathbb{P}^{2p})^m$ computed component-wise, $m$-vectors $\pmb{\omega}_i$ of nonlinear coefficients in each cell computed component-wise, while of course the linear coefficients $\mu^i_{j,\alpha}$ are independent of the component being reconstructed.

\subsection{The space-time first order implicit prediction} \label{sec:quinpi:predictor}

Without loss of generality, assume that the nodes $c_1,\dots,c_s$ of the DIRK method are ordered. Then, we approximate the system of ODEs~\eqref{eq:spaceApprox} with an implicit first order scheme at any time $t^{(k)}=(n+c_{k}) \Delta t \in [n\Delta t,(n+1)\Delta t]$, where $\Delta t$ is the time-step of the high order scheme~\eqref{eq:implicit}. Specifically, the system is numerically approximated in space using piecewise constant reconstructions and integrated in time using a composite backward Euler method, providing the $s$ approximations 
$\ca{\bigvec{U}}^{\star,(k)} = \left[\ca{\mathbf{U}}^{\star,(k)}_1,\dots,\ca{\mathbf{U}}^{\star,(k)}_N\right]^T$, $k=1,\dots,s$. Therefore, for each $j=1,\dots,N$, 
$\ca{\mathbf{U}}_j^{\star,(k)} \approx \ca{\mathbf{u}}_j(t^{(k)})$ is given by
\begin{subequations} \label{eq:predictor}
	\begin{align}
		\ca{\mathbf{U}}_j^{\star,(k)} &= \ca{\mathbf{U}}_j^{\star,(k-1)} - \frac{\theta_k\Delta t}{h} \left[ \NumFlux_{j+\frac12}^{\star,(k)} - \NumFlux_{j-\frac12}^{\star,(k)} \right], \quad k=1,\dots,s, \label{eq:predictor:stage} \\
		\ca{\mathbf{U}}_j^{\star,n+1} &= \ca{\mathbf{U}}_j^{n} - \frac{\Delta t}{h} \sum_{k=1}^s \theta_k \left[ \NumFlux_{j+\frac12}^{\star,(k)} - \NumFlux_{j-\frac12}^{\star,(k)} \right], \quad n\geq 0, \label{eq:predictor:update} \\
		\NumFlux_{j+\frac12}^{\star,(k)} &= \NumFlux\left(\ca{\mathbf{U}}_{j}^{\star,(k)},\ca{\mathbf{U}}_{j+1}^{\star,(k)}\right) \in \R^m \label{eq:predictor:flux}
	\end{align}
\end{subequations}
where $\theta_k:=c_{k}-c_{k-1}$, with $c_0=0$ and $\ca{\mathbf{U}}_j^{\star,(0)}:=\ca{\mathbf{U}}_j^{n}$, and where the convention 
$$
\ca{\mathbf{U}}_{0}^{\star,(k)} = u_{\mathsf{out}}^-, \quad \ca{\mathbf{U}}_{N+1}^{\star,(k)} = u_{\mathsf{out}}^+,
$$
is used. The outer values $u_{\mathsf{out}}^{\mp}$ are imposed by the chosen boundary conditions. For instance, for periodic conditions one has
\begin{equation} \label{eq:predictor:periodic}
	\ca{\mathbf{U}}_{0}^{\star,(k)} = \ca{\mathbf{U}}_{N}^{\star,(k)}, \quad \ca{\mathbf{U}}_{N+1}^{\star,(k)} = \ca{\mathbf{U}}_{1}^{\star,(k)},
\end{equation}
whereas for free-flow conditions one gets
\begin{equation} \label{eq:predictor:freeflow}
	\ca{\mathbf{U}}_{0}^{\star,(k)} = \ca{\mathbf{U}}_{1}^{\star,(k)}, \quad \ca{\mathbf{U}}_{N+1}^{\star,(k)} = \ca{\mathbf{U}}_{N}^{\star,(k)}.
\end{equation}
For more general boundary conditions, one should project on characteristic variables, providing the prescribed upwind conditions.
Summarizing, the first order scheme~\eqref{eq:predictor} is equivalent to applying a DIRK scheme with Butcher tableau given by
\begin{equation}
	\label{eq:tableau:be}
	\begin{array}{c|cccc}
		c_1 & \theta_{1} & 0 & \dots & 0 \\[1.5ex]
		c_2 & \theta_{1} & \theta_{2} & \dots & 0 \\[1.5ex]
		\vdots & \vdots & \vdots & \ddots & \\[1.5ex]
		c_s & \theta_{1} & \theta_{2} & \dots & \theta_{s} \\[1ex]
		\hline
		&&&&\\[-1.8ex]
		& \theta_1 & \theta_2 & \dots & \theta_s
	\end{array}
\end{equation}
where the coefficients $c_i$, $i=1,\dots,s$, are the same as the high order DIRK in~\eqref{eq:tableau:dirk}.

Notice that the numerical flux function $\NumFlux$ in~\eqref{eq:predictor:flux} is now computed on piecewise constant, unlimited, reconstructions from the cell averages. In fact, first order schemes do not require space-limiting, because they are unconditionally Total Variation Diminishing, see~\cite[Section 2]{PSV23:Quinpi}. Therefore, despite of the high order approximation~\eqref{eq:implicit}, the first order scheme~\eqref{eq:predictor} is characterized by a single nonlinearity, that is the one induced by the flux function $\mathbf{f}$.
Computing~\eqref{eq:predictor:stage} requires the solution of $s$ nonlinear systems of dimension $m N$, which are fully linear with respect to the unknown cell averages:
\begin{equation} \label{eq:predictor:system}
	\mathbf{G}^\star(\ca{\bigvec{U}}^{\star,(k)}) := 
	\ca{\bigvec{U}}^{\star,(k)} 
	+ \frac{\theta_k \Delta t}{h} \Delta\bigvec{F}^{\star,(k)} 
	- \ca{\bigvec{U}}^{\star,(k-1)} = \mathbf{0},
\end{equation}
where $\Delta\bigvec{F}^{\star,(k)}\in\R^{mN}$ whose $j$-th block is $\Delta\bigvec{F}^{\star,(k)}_j=\left(\NumFlux_{j+\nicefrac12}^{\star,(k)}-\NumFlux_{j-\nicefrac12}^{\star,(k)} \right) \in\R^m$, for $j=1,\dots,N$.

The solution of~\eqref{eq:predictor:system} requires the use of a nonlinear solver $s$ times within a single time-step. In this work, we rely on Newton's iterations
\begin{equation} \label{eq:predicor:newton}
	\ca{\bigvec{U}}^{\star,(k)}_{(\ell+1)} 
	= \ca{\bigvec{U}}^{\star,(k)}_{(\ell)} 
	- \left( \mathbb{I} + \frac{\theta_k \Delta t}{h} J^{\star}\left( \ca{\bigvec{U}}^{\star,(k)}_{(\ell)} \right) \right)^{-1} \mathbf{G}^\star\left( \ca{\bigvec{U}}^{\star,(k)}_{(\ell)} \right), \quad \ell \geq 0,
\end{equation}
with given initial guess $\ca{\bigvec{U}}^{\star,(k)}_{(0)}=\ca{\bigvec{U}}^{\star,(k-1)}$. Here, $\mathbb{I}\in\R^{mN\times mN}$ is the identity matrix, whereas $J^{\star}$ is the Jacobian matrix of the numerical flux difference $\Delta \bigvec{F}^{\star}$ which is organized in $N\times N$ blocks of size $m\times m$. Then, each block of the Jacobian is given by
$$
\left( J^{\star} \right)_{ji} = \frac{\partial \Delta \bigvec{F}^{\star}_j}{\partial \ca{\mathbf{U}}^{\star}_i} \in \R^{m\times m}, \quad j,i=1,\dots,N.
$$
In particular, for a first order space approximation one has
%\begin{align*}
%	\Delta \NumFlux^{\star}_1 &= \Delta \NumFlux^{\star} (\ca{\mathbf{U}}^{\star}_{1},\ca{\mathbf{U}}^{\star}_2,\ca{\mathbf{U}}^{\star}_{N}), \\
%	\Delta \NumFlux^{\star}_j &= \Delta \NumFlux^{\star} (\ca{\mathbf{U}}^{\star}_{j-1},\ca{\mathbf{U}}^{\star}_j,\ca{\mathbf{U}}^{\star}_{j+1}), \quad j=2,\dots,N-1 \\
%	\Delta \NumFlux^{\star}_N &= \Delta \NumFlux^{\star} (\ca{\mathbf{U}}^{\star}_{1},\ca{\mathbf{U}}^{\star}_{N-1},\ca{\mathbf{U}}^{\star}_{N}),
%\end{align*}
%which are fully linear on the data, it turns out
that the Jacobian is a tridiagonal block matrix of the form
$$
\frac{\partial \Delta \bigvec{F}^{\star}_j}{\partial \ca{\mathbf{U}}^{\star}_i} = \begin{cases}
	-\dfrac{\partial \NumFlux}{\partial \mathbf{v}} \left( \ca{\mathbf{U}}^{\star}_i,\ca{\mathbf{U}}^{\star}_{i+1} \right), & \text{if } i=j-1 \\[2.5ex]
	\dfrac{\partial \NumFlux}{\partial \mathbf{v}}\left( \ca{\mathbf{U}}^{\star}_i,\ca{\mathbf{U}}^{\star}_{i+1} \right) - \dfrac{\partial \NumFlux}{\partial \mathbf{w}} \left( \ca{\mathbf{U}}^{\star}_{i-1},\ca{\mathbf{U}}^{\star}_i \right), & \text{if } i=j \\[2.5ex]
	\dfrac{\partial \NumFlux}{\partial \mathbf{w}} \left( \ca{\mathbf{U}}^{\star}_{i-1},\ca{\mathbf{U}}^{\star}_i \right), & \text{if } i=j+1 \\
\end{cases}
$$
for $j=2,\dots,N-1$, where $\partial_{\mathbf{v}} \NumFlux$ and $\partial_{\mathbf{w}} \NumFlux$ denote the Jacobian of the vector-valued numerical flux function~\eqref{eq:numfluxfnc} with respect to its first and second variable, respectively.

\subsection{The space-time third order implicit correction} \label{sec:quinpi:corrector}

Assume that the first order predictor is available at a time $t^{(k)} \in [n\Delta t,(n+1)\Delta t]$. The obtained approximation $\ca{\bigvec{U}}^{\star,(k)}$ is used to evaluate the nonlinear terms of the high order BED as explained at the beginning of this section, see~\eqref{eq:bed:linearized}. Therefore, the third order ``correction'' at time $t^{(k)}$ can be computed as in~\eqref{eq:implicit:stage} with
\begin{equation} \label{eq:implicit:fluxLinearized}
	\begin{aligned}
		\NumFlux_{j+\frac12}^{(k)} = \NumFlux\left(\mathbf{U}_{j+\frac12}^{-,(k)},\mathbf{U}_{j+\frac12}^{+,(k)}\right) \approx \hat{\NumFlux}_{j+\frac12}^{(k)} = \NumFlux\left(\hat{\mathbf{U}}_{j+\frac12}^{-,(k)},\hat{\mathbf{U}}_{j+\frac12}^{+,(k)}\right),
	\end{aligned}
\end{equation}
where $\hat{\mathbf{U}}_{j+\frac12}^{\mp,(k)}$ are the linearized BED, which are fully linear with respect to the unknown cell averages at time $t^{(k)}=t_n+c_k\Delta t$. Therefore, the nonlinear system
\begin{equation} \label{eq:corrector:system}
	\mathbf{G}\left(\ca{\bigvec{U}}^{(k)}\right) 
	:= \ca{\bigvec{U}}^{(k)} 
	+ \frac{a_{kk}\Delta t}{h} \Delta \hat{\bigvec{F}}^{(k)} 
	- \ca{\bigvec{U}}^{(k-1)} = \mathbf{0},
\end{equation}
obtained from~\eqref{eq:implicit:stage} and \eqref{eq:implicit:fluxLinearized}, can be tackled with Newton's iterations
\begin{equation} \label{eq:corrector:newton}
	\ca{\bigvec{U}}^{(k)}_{(\ell+1)} = \ca{\bigvec{U}}^{(k)}_{(\ell)} 
	- \left( \mathbb{I} + \frac{a_{kk} \Delta t}{h} J\left( \ca{\bigvec{U}}^{(k)}_{(\ell)} \right) \right)^{-1} \mathbf{G}\left( \ca{\bigvec{U}}^{(k)}_{(\ell)} \right), \quad \ell \geq 0,
\end{equation}
with initial guess $\ca{\bigvec{U}}^{(k)}_{(0)}$ set to the predictor's output. The matrix $J\in\R^{mN\times mN}$ is the Jacobian of the system~\eqref{eq:corrector:system} and one has that each block $\left( J \right)_{ji} \in \R^{m\times m}$, $j,i=1,\dots,N$, can be written as
\begin{align*}
	\left( J \right)_{ji} = \frac{\partial \Delta \hat{\bigvec{F}}_j}{\partial \ca{\mathbf{U}}_{i}}
	=& \frac{\partial \NumFlux}{\partial \mathbf{v}} \left( \hat{\mathbf{U}}^-_{j+\frac12} , \hat{\mathbf{U}}^+_{j+\frac12}\right) \frac{\partial \hat{\mathbf{U}}^-_{j+\frac12}}{\partial \ca{\mathbf{U}}_i} + \frac{\partial \NumFlux}{\partial \mathbf{w}} \left( \hat{\mathbf{U}}^-_{j+\frac12} , \hat{\mathbf{U}}^+_{j+\frac12}\right) \frac{\partial \hat{\mathbf{U}}^+_{j+\frac12}}{\partial \ca{\mathbf{U}}_i} \\
	&- \frac{\partial \NumFlux}{\partial \mathbf{v}} \left( \hat{\mathbf{U}}^-_{j-\frac12} , \hat{\mathbf{U}}^+_{j-\frac12}\right) \frac{\partial \hat{\mathbf{U}}^-_{j-\frac12}}{\partial \ca{\mathbf{U}}_i} - \frac{\partial \NumFlux}{\partial \mathbf{w}} \left( \hat{\mathbf{U}}^-_{j-\frac12} , \hat{\mathbf{U}}^+_{j-\frac12}\right) \frac{\partial \hat{\mathbf{U}}^+_{j-\frac12}}{\partial \ca{\mathbf{U}}_i}
\end{align*}
where $\partial_{\mathbf{v}}\NumFlux$ and $\partial_{\mathbf{w}}\NumFlux$ contain the Jacobian of the flux function $\mathbf{f}$, whereas $\partial_{\ca{\mathbf{U}}_i} \mathbf{U}^{\mp}_{j\pm\nicefrac12}$ are assembled with the frozen nonlinear weights, cf.~\eqref{eq:bed:linearized}. Observe that the width of the band of the block-banded matrix that defines the Jacobian $J$ is still as in the fully implicit approach, because of the stencil of the space reconstructions, but the entries can now be computed explicitly.

For nontrivial boundary conditions, the above expression for the Jacobian blocks has to be modified according to the dependence of the external BED ($\ca{\mathbf{U}}_{\nicefrac12}^-$ and $\ca{\mathbf{U}}_{N+\nicefrac12}^+$) on the interior data.

Finally, once all the stage values $\ca{\bigvec{U}}^{(k)}$, $k=1,\dots,s$, are computed, the third order solution $\ca{\bigvec{U}}^{n+1}$ at time level $(n+1)\Delta t$ is obtained as in~\eqref{eq:implicit:update}.

\begin{remark}[Implicit reconstruction along characteristic variables.]
	The Quinpi framework can be used also for reconstruction along characteristic variables.
	The linearity of the scheme can be kept provided that the nonlinear transformation matrix from conservative to characteristic variables is computed by the predictor scheme and stored at each time-step $t^{(k)} = t^{n} + c_k \Delta t$, $k=1,\dots,s$.
	
	More precisely, let $\mathbf{U}$ denote the vector of conservative variables and let $\mathbf{V}$ denote the vector of characteristic variables which are linked through the nonlinear transformation
	$$
	\mathbf{V} = T(\mathbf{U}).
	$$
	Then, the time advancement $\mathbf{\ca{U}}_j^{(k-1)} \mapsto \mathbf{\ca{U}}_j^{(k)}$ on the cell $\Omega_j$, with reconstruction from cell averages performed in characteristic variables, is obtained for $k=2,\dots,s$ as follows.
	
	First one can perform the time advancement of the predictor solution $\mathbf{\ca{U}}_j^\star$ and store the nonlinear transformation matrix:
	$$
	\{ \mathbf{\ca{U}}_j^{\star,(k-1)} \} \xmapsto{\textnormal{advance}} \{ \mathbf{\ca{U}}_j^{\star,(k)} \} \xmapsto{\textnormal{store}} \{ T_j^{(k)} \},
	$$
	where $T_j^{(k)}$ is the map to characteristic variables on the state $\mathbf{\ca{U}}^{\star,(k)}_j$, i.e.~the map whose columns are the right eigenvectors of the flux Jacobian at $\mathbf{\ca{U}}^{\star,(k)}_j$.
	
	Then, the high-order solution at time $t^{(k)}$ is obtained as
	$$
	\{ \mathbf{\ca{U}}_\ell^{(k-1)},\mathbf{\ca{U}}_\ell^{\star,(k)} \}_{\ell\in\mathcal{S}_{\text{opt}}}
	\xmapsto{T_j^{(k)}, \textnormal{ linear.~BED}} 
	\{ \mathbf{\hat{V}}_{j\pm\frac12}^{(k)} \} \xmapsto{(T_j^{(k)})^{-1}} \{ \mathbf{\hat{U}}_{j\pm\frac12}^{(k)} \} \xmapsto{\textnormal{advance}} \{ \mathbf{\ca{U}}_j^{(k)} \}.
	$$
\end{remark}

\subsubsection{On the accuracy of the reconstruction in the smooth case}

In CWENO reconstructions optimal accuracy is obtained provided that the nonlinear weights $\omega_k$ are sufficiently close to the linear weights $d_k$. More precisely, for the third order CWENO scheme a sufficient condition to have optimal accuracy on smooth data is, cf.~\cite{CPSV:cweno},
$$
d_k-\omega_k = \mathcal{O}(h), \quad k=0,L,R.
$$
Indeed, let $R_j^{\text{CWZ}}\in\mathbb{P}^2$ be the reconstruction polynomial, defined in~\eqref{eq:precCWZ}, of a sufficiently smooth function $u$. Then, the reconstruction error is
\begin{equation} \label{eq:rec:accuracy}
	\begin{aligned}
		u(x,\cdot)-R_j^{\text{CWZ}}(x;\cdot) =& \underbrace{(u(x,\cdot)-\Popt(x;\cdot))}_{\mathcal{O}(h^{3})} \\
		&+ \sum_{k=L,R} (d_k-\omega_k) \underbrace{(P_k(x;\cdot)-u(x,\cdot))}_{\mathcal{O}(h^{2})}, \quad x\in\Omega_j,
	\end{aligned}
\end{equation}
when the polynomials $\Popt$, $P_L$ and $P_R$ interpolate the exact cell averages of $u$. Therefore, the optimal order of accuracy, in this case $3$, is obtained on smooth data if
\begin{equation} \label{eq:wmend}
	d_k-\omega_k = \mathcal{O}(h), \quad k=0,L,R.
\end{equation}

We investigate numerically the convergence rate $d_k-\omega_k$ when the nonlinear weights are computed on a low order approximation of the cell-averages of a smooth function $u$. Thus, let us consider
$$
u(x) = \sin(\pi x)+\sin(15\pi x) \exp(-20x^2), \quad x\in[-1,1].
$$
We discretize the space domain by $N$ cells with periodic boundary conditions. We take the optimal weights $d_0 = 0.75$, $d_L = d_R = 0.125$. Then, we compute the nonlinear weights $\hat{\omega}_{j,k}$, $k=0,L,R$, based on a first order approximation of the cell averages of $u$ in the $j$-th cell, and the nonlinear weights $\omega_{j,k}$, $k=0,1,2$, on a third order approximation of the cell averages of $u$ in the $j$-th cell. We study the behavior of the errors
\begin{equation*} \label{eq:err3}
	\begin{aligned}
		\widehat{\text{err}}(N) &= \frac{1}{N} \sum_{j=1}^N \max_{k=0,L,R} |d_k-\hat{\omega}_{j,k}| \\
		\text{err}(N) &= \frac{1}{N} \sum_{j=1}^N \max_{k=0,L,R} |d_k-\omega_{j,k}|
	\end{aligned}
\end{equation*}
under grid refinement. The results are given in Table~\ref{tab:order3_omegad} which shows that the sufficient condition on the convergence rate of the nonlinear weights to the linear weights to guarantee optimal accuracy is fulfilled in both cases. In Quinpi the accuracy conditions in~\eqref{eq:rec:accuracy} are guaranteed since the $\Popt$ and $P_k$ polynomials are computed on the high order cell averages in the nonlinear solver, while the first order predictor is accurate enough to guarantee condition~\eqref{eq:wmend} on the nonlinear weights.

\begin{table}[t!]
	\centering
	\caption{Experimental convergence of the nonlinear to the linear weights.\label{tab:order3_omegad}}
	\begin{tabular}{ccccc}
		\toprule
		$N$ & $\widehat{\text{err}}$ & rate & $\text{err}$ & rate \\
		\midrule
		20 & 3.58$e-1$ & -      & 3.83$e-1$ & - \\
		40 & 2.56$e-1$ & 0.48   & 2.92$e-1$ & 0.39 \\
		80 & 2.43$e-1$ & 0.08   & 2.46$e-1$ & 0.25 \\
		160 & 1.57$e-1$ & 0.63  & 1.56$e-1$ & 0.65 \\
		320 & 7.23$e-2$ & 1.12  & 7.23$e-2$ & 1.11 \\
		640 & 3.96$e-2$ & 0.87  & 3.96$e-2$ & 0.87 \\
		1,280 & 1.76$e-2$ & 1.17 & 1.76$e-2$ & 1.17 \\
		2,560 & 5.59$e-3$ & 1.65 & 5.59$e-3$ & 1.65 \\
		5,120 & 1.37$e-3$ & 2.03 & 1.37$e-3$ & 2.03 \\
		\bottomrule
	\end{tabular}
\end{table}

\subsection{Conservative a-posteriori time-limiting based on numerical entropy} \label{sec:quinpi:entropy}

The solution obtained with the predictor-corrector approach is third order accurate, and has some control over spurious oscillations thanks to the limited CWENO space reconstruction. However, when using a large $\Delta t$, space limiting is not enough, and the high order solution may still exhibit oscillations as already noticed in~\cite{PSV23:Quinpi,2020Arbogast,2007DurasaisamyBaeder}. In such situations, limiting in time is also required.
For this reason, the Quinpi scheme that we propose in this paper employs a time-limiting procedure of the third order solution.

In~\cite{PSV23:Quinpi} the time-limiting is performed in a WENO-like framework. In fact, nonlinear weights are suitably defined in order to blend cell averages between the computed high order solution and the low order predictor, which is reliable, stable and non-oscillatory. However, the procedure in~\cite{PSV23:Quinpi} does not preserve mass conservation property and must be followed by a suitable conservative correction, analogously to~\cite{Marsha:AMR} in adaptive mesh refinement.

Extending the conservative correction of~\cite{PSV23:Quinpi} to a generic system~\eqref{eq:hyp:sys} is nontrivial. Instead, conservation can be ensured avoiding a cell-centered blending by means of flux-based Runge-Kutta, as, e.g., in~\cite{Ketcheson:fluxbased:2013,2020Arbogast}.
Similarly, a flux-centered conservative a-posteriori time-limiting inspired by the MOOD technique~\cite{CDL11:MOOD,CDL12:MOOD} was investigated in~\cite{VTSP23:Quinpi:Book}. MOOD was originally designed as an a-posteriori space-limiting technique for multi-dimensional finite volume schemes. Instead, in~\cite{VTSP23:Quinpi:Book}, the typical MOOD detectors were used to limit the high order solution at time level $(n+1)\Delta t$. However, contrary to the classical MOOD approach, in~\cite{VTSP23:Quinpi:Book} the method uses a convex combination of the high order and the low-order numerical fluxes at the interfaces of oscillatory cells. A similar idea was also employed in~\cite{EBS22:Implicit:Networks} for a-posteriori limiting of fully implicit finite volume schemes on transport networks.

In this paper we still rely on the MOOD technique, and the high order numerical fluxes of oscillatory cells are dropped to the low-order numerical fluxes of the predictor. The novelty here is that the detection of troubled cells is performed with the numerical entropy production introduced in~\cite{PS11:numerical:entropy,Puppo04:numerical:entropy}, extended to balance laws in~\cite{PS16:entropy:balance}, and already exploited in \cite{SCR:CWENOquadtree,SL:18:AMRMOOD} as a-posteriori error of adaptive schemes and indicator of the qualitative structure of the flow.

\subsubsection{Indicators based on numerical entropy production} \label{sec:quinpi:entropy:indicators}

Assume that system~\eqref{eq:hyp:sys} is endowed with an entropy-entropy flux pair, namely that there exists a convex function $\mathbf{u}\in\R^m\mapsto\eta(\mathbf{u})\in\R$ and a corresponding entropy flux $\mathbf{u}\in\R^m\mapsto\psi(\mathbf{u})\in\R$ such that $\nabla^T\eta(\mathbf{u}) \mathbf{f}^\prime(\mathbf{u})=\nabla^T\psi(\mathbf{u})$. Then, all entropy solutions of~\eqref{eq:hyp:sys} satisfy the entropy inequality
\begin{equation} \label{eq:entropy}
	\frac{\partial}{\partial t} \eta(\mathbf{u}(x,t)) + \frac{\partial}{\partial x} \psi(\mathbf{u}(x,t)) \leq 0,
\end{equation}
in the weak sense. It is well known that~\eqref{eq:entropy} is an equality on smooth flows. Instead, if the solution has a singularity, the sign of~\eqref{eq:entropy} selects the unique physically admissible weak solution of~\eqref{eq:hyp:sys}. For this reason, the novel idea of~\cite{Puppo04:numerical:entropy} was to consider the numerical residual of
the entropy inequality as an a-posteriori error indicator for the numerical solution of~\eqref{eq:hyp:sys}, since it is a scalar value (even in the case of systems of conservation laws) which provides information on the size of the local truncation error and on the presence of singularities. Recently, the numerical entropy residual has been employed also in the context of stochastic Galerkin formulations of hyperbolic systems, see~\cite{2023GersterSemplice}.

Following~\cite{PS11:numerical:entropy} we give the subsequent definition.

\begin{definition}
	Let $\ca{\mathbf{U}}^n$ be the solution of~\eqref{eq:spaceApprox} obtained using a $s$-stage Runge-Kutta scheme of order $p$ with weights $\{ b_i \}_{i=1}^s$ and time-step $\Delta t$. Then, the numerical entropy production $\left(S^p\right)_j^n$ of the scheme in the control volume $V_j^n = \Omega_j \times [n\Delta t , (n+1)\Delta t]$ is given by
	\begin{equation} \label{eq:numerical:entropy}
		\left(S^p\right)_j^n = \frac{1}{\Delta t} \left[ \ca{Q}\left( \eta\left( \ca{\mathbf{U}}^{n+1} \right) \right)_j - \ca{Q}\left( \eta\left( \ca{\mathbf{U}}^{n} \right) \right)_j + \frac{\Delta t}{h} \sum_{i=1}^s b_i \left( \Psi_{j+\frac12}^{(i)} - \Psi_{j-\frac12}^{(i)} \right) \right]
	\end{equation}
	where
	$$
	\Psi_{j+\frac12}^{(i)} = \Psi\left( \mathbf{U}_{j+\frac12}^{-,(i)},\mathbf{U}_{j+\frac12}^{+,(i)} \right),
	$$
	is a numerical entropy flux function consistent with the exact entropy flux $\psi$ of~\eqref{eq:entropy}, with $(\mathbf{v},\mathbf{w}) \in \R^m\times\R^m \to \Psi(\mathbf{v},\mathbf{w})\in\R$, and where $\ca{\mathbf{U}}_{j\pm\nicefrac12}^{\mp,(i)}$ are the BED of the numerical solution. In~\eqref{eq:numerical:entropy}, $\ca{Q}\left( \cdot \right)_j$ denotes a quadrature rule of order $p$ to compute the cell average on the cell $\Omega_j$.
\end{definition}

The rate of convergence of~\eqref{eq:numerical:entropy} is shown in~\cite{PS11:numerical:entropy}. In particular, for a general numerical scheme of order $p$ one has
$$
\left(S^p\right)_j^n = \begin{cases}
	\mathcal{O}(h^p), & \text{if the solution is regular on $\Omega_j$} \\
	\mathcal{O}(h^{-1}), & \text{if the solution has a shock on $\Omega_j$} \\
	\mathcal{O}(h), & \text{if $\Omega_j$ exhibits a rarefaction corner} \\
	\mathcal{O}(1), & \text{if $\Omega_j$ exhibits a contact wave}.
\end{cases}
$$
In other words, $\left(S^p\right)_j^n\to 0$, as $h\to 0$, with the same rate of the local truncation error of the scheme on smooth flows. Whereas, the numerical entropy production increases in magnitude on shocks.

For the first order predictor solution $\ca{\bigvec{U}}^\star$, defined in Section~\ref{sec:quinpi:predictor}, we can consider the midpoint quadrature rule and observe that the cell average and the pointwise value at the cell center are $\mathcal{O}(h^2)$ apart, and thus set
$
\ca{Q}\left( \eta\left( \ca{\mathbf{U}}^{n} \right) \right)_j =  \eta\left(\ca{\mathbf{U}}_j^{n}\right),
$
in~\eqref{eq:numerical:entropy}. The numerical entropy production becomes
\begin{equation} \label{eq:entropy:1}
	\left(S^1\right)_j^n = \frac{1}{\Delta t} \left[ \eta\left( \ca{\mathbf{U}}^{\star,n+1}_j \right) - \eta\left( \ca{\mathbf{U}}^{\star,n}_j \right) + \frac{\Delta t}{h} \sum_{i=1}^s \theta_i \left( \Psi_{j+\frac12}^{\star,(i)} - \Psi_{j-\frac12}^{\star,(i)} \right) \right],
\end{equation}
where
$$
\Psi_{j+\frac12}^{\star,(i)} = \Psi\left( \mathbf{U}_j^{\star,(i)},\mathbf{U}_{j+1}^{\star,(i)} \right).
$$
Instead, for the third order solution $\ca{\bigvec{U}}$ defined in Section~\ref{sec:quinpi:corrector}, in this work we consider the two-point Gauss-Legendre quadrature rule:
$$
\ca{Q}\left( \eta\left( \ca{\mathbf{U}}^{n} \right) \right)_j = \frac{1}{2} \left( \eta\left(\mathbf{U}_{j-\frac{\sqrt{3}}{6}}^{n}\right) + \eta\left(\mathbf{U}_{j+\frac{\sqrt{3}}{6}}^{n}\right) \right),
$$
where $\mathbf{U}_{j\pm\nicefrac{\sqrt{3}}{6}}^{n}$ denote the reconstructions of the solution in $x_j\pm\nicefrac{\sqrt{3}h}{6}$. Observe that the computation of the quadrature rule does not require significant additional computational cost since the computation of $\mathbf{U}_{j\pm\nicefrac{\sqrt{3}}{6}}^{n}$ involves just the evaluation of a CWENO reconstruction polynomial, which is already available, because it is computed in order to advance the solution from $n\Delta t$ to $(n+1)\Delta t$. Thus, for the third order solution, the numerical entropy production becomes
\begin{equation} \label{eq:entropy:3}
	\begin{aligned}
		\left(S^3\right)_j^n = \frac{1}{\Delta t}  
		\bigg[ &\frac{1}{2} \left( \eta\left(\mathbf{U}_{j-\frac{\sqrt{3}}{6}}^{n+1}\right) + \eta\left(\mathbf{U}_{j+\frac{\sqrt{3}}{6}}^{n+1}\right) \right) - \frac{1}{2} \left( \eta\left(\mathbf{U}_{j-\frac{\sqrt{3}}{6}}^{n}\right) + \eta\left(\mathbf{U}_{j+\frac{\sqrt{3}}{6}}^{n}\right) \right)  \\
		& + \frac{\Delta t}{h} \sum_{i=1}^s b_i \left( \Psi_{j+\frac12}^{(i)} - \Psi_{j-\frac12}^{(i)} \right) \bigg].
	\end{aligned}
\end{equation}

At the discrete level, $\left(S^1\right)_j^n$ and $\left(S^3\right)_j^n$ provide a regularity indicator of the predictor and of the corrector solution, respectively. Therefore, in order to determine where time-limiting has to occur, we detect spurious oscillations using $\left(S^1\right)_j^n$ and $\left(S^3\right)_j^n$ in the following ways.
\begin{description}
	\item[$\mathcal{I}_1$:] The indicator~\eqref{eq:entropy:3} of the third order scheme is used to detect troubled cells. We expect that on smooth cells $\Omega_j$ one has $\left(S^3\right)_j^n = \mathcal{O}(h^3)$, whereas $\left(S^3\right)_j^n = \mathcal{O}(h^{-1})$ in presence of shocks. In particular, a cell $\Omega_j$ is marked for time-limiting if
	\begin{equation} \label{eq:strategy1}
		\left(S^3\right)_j^n > \gamma_{1},
	\end{equation}
	where $\gamma_{1}$ is a given threshold.
	\item[$\mathcal{I}_2$:] The ratio of the indicators~\eqref{eq:entropy:1} and~\eqref{eq:entropy:3} is used to detect problematic cell $\Omega_j$ whenever
	\begin{equation} \label{eq:strategy2}
		\frac{\left(S^3\right)_j^n}{\left(S^1\right)_j^n+\sigma} > \gamma_{2},
	\end{equation}
	where $\gamma_{2}$ is a given threshold and $\sigma$ is a small quantity which prevents $\nicefrac{\left(S^3\right)_j^n}{\left(S^1\right)_j^n} \sim \mathcal{O}(1)$ on constant regions of the solution. In fact, we expect that on a smooth cell $\Omega_j$ one has $\nicefrac{\left(S^3\right)_j^n}{\left(S^1\right)_j^n} = \mathcal{O}(h^2)$, whereas $\nicefrac{\left(S^3\right)_j^n}{\left(S^1\right)_j^n} = \mathcal{O}(1)$ if the solution is not regular on $\Omega_j$.
	\item[$\mathcal{I}_3$:] Both the detecting techniques $\mathcal{I}_1$ and $\mathcal{I}_2$ are used, so that a cell $\Omega_j$ is marked for time-limiting if
	\begin{equation} \label{eq:strategy3}
		\left(S^3\right)_j^n > \gamma_{1} \ \text{ and } \ \frac{\left(S^3\right)_j^n}{\left(S^1\right)_j^n+\sigma} > \gamma_{2}.
	\end{equation}
	We will compare these strategies in the numerical section. For now let us remark that $\mathcal{I}_2$ should be more robust than $\mathcal{I}_1$ since it involves ratios, but that it may fail on flat regions where both $S^3$ and $S^1$ can be floating point zeros. Thus in $\mathcal{I}_3$ we propose to combine the $\mathcal{I}_1$ and the $\mathcal{I}_2$ strategies. In fact, the $\mathcal{I}_2$ strategy may mark flat regions as irregular cells, i.e.~where a shock appears, if $\left(S^3\right)_j^n = \mathcal{O}(\left(S^1\right)_j^n)$. Then, the $\mathcal{I}_1$ strategy allows to detect the regularity of the solution. 
\end{description}

The thresholds $\gamma_1$ and $\gamma_2$ are very important for the performance of the regularity detection. We will come back to this point in Section~\ref{sec:numerics} when we discuss the values we choose for the two thresholds.

\subsubsection{Time limited solution and the final Quinpi algorithm}

\begin{figure}[t!]
	\centering
	\includegraphics[width=\textwidth]{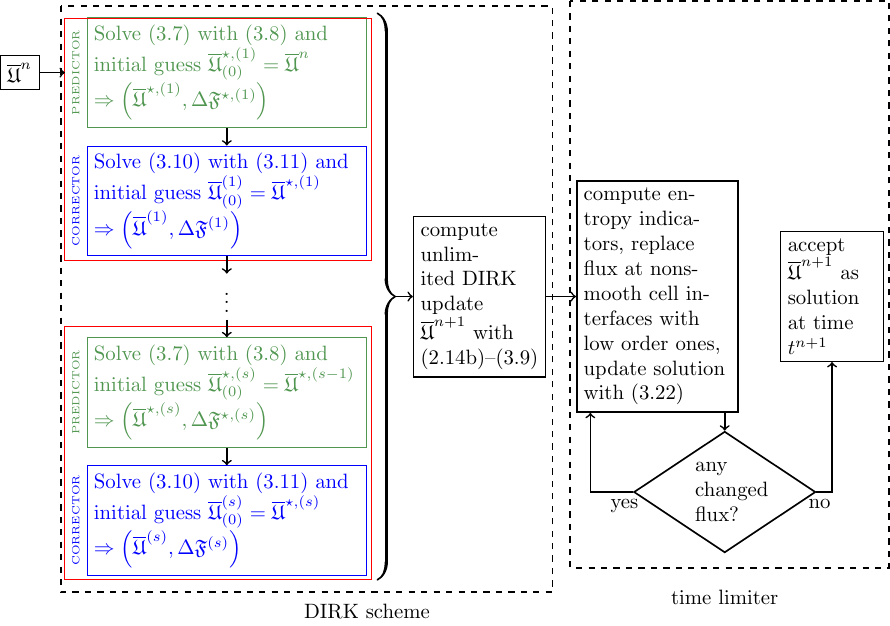}
	\caption{The algorithm of the Quinpi scheme for the computation of the solution in a single time-step. The DIRK routine, which provides the high order solution with the space limiters computed by the predictor, is followed by the time limiter which downgrades the numerical fluxes at the interface of an irregular cell. \label{fig:algorithm}}
\end{figure}

Once the detector has been chosen, the flux-centered time-limiting is performed by computing the limited numerical flux $\NumFlux^{\mathsf{TL},(k)}_{j+\nicefrac12}$
for the face between $\Omega_{j}$ and $\Omega_{j+1}$, for 
all stages $k=1,\dots,s$, as
\begin{equation} \label{eq:timelimited:fluxes}
	\NumFlux^{\mathsf{TL},(k)}_{j+\frac12} = \begin{cases}
		\theta_k \NumFlux^{\star,(k)}_{j+\frac12}, & \text{if either $\Omega_j$ or $\Omega_{j+1}$ is marked},\\[2ex]
		b_k \NumFlux^{(k)}_{j+\nicefrac12}, & \text{if both $\Omega_j$ and $\Omega_{j+1}$ are not marked},
	\end{cases}
\end{equation}
where the acronym $\mathsf{TL}$ stands for Time-Limited. Notice that, since the predictor scheme is based on a composite backward Euler, the low order numerical fluxes $\NumFlux^{\star,(k)}_{j\pm\nicefrac12}$ are available at each stage $k=1,\dots,s$, for $j=1,\dots,N$.

Thus, the high order numerical fluxes at the interfaces of a problematic cell $\Omega_j$ are replaced with the low order numerical fluxes reducing locally the order of the solution, which is then computed as
\begin{equation} \label{eq:timelimited:sol}
	\ca{\mathbf{U}}_j^{n+1} = \ca{\mathbf{U}}_j^{n} - \frac{\Delta t}{h} \sum_{k=1}^s \left[ \NumFlux_{j+\tfrac12}^{\mathsf{TL},(k)} - \NumFlux_{j-\tfrac12}^{\mathsf{TL},(k)} \right], \quad j=1,\dots,N.
\end{equation}
The time-limiting procedure is repeated until the solution $\ca{\bigvec{U}}^{n+1}$ is detected as smooth on each cell.

Overall, the update of the solution from time $n\Delta t$ to time $(n+1)\Delta t$ computed with the Quinpi scheme is summarized in Figure~\ref{fig:algorithm}. 

\section{Numerical simulations} \label{sec:numerics}

All the schemes we consider in this section use the third order CWENOZ reconstruction without ghost cells, see~\cite{STP23:cweno:boundary}, described in Definition~\ref{def:CWENOZb} and \ref{def:CWENOZAO} with parameters $q=p=2$ in \eqref{eq:omegaZ} and
\begin{equation*}
	\begin{cases}
		d_L=\frac18, \ d_R=\frac18, & \text{for the inner computational cells} \\
		\tilde{d}=\max\{h,0.01\}, \ d=\frac14, & \text{for the first/last computational cell.}
	\end{cases}
\end{equation*}

For the purposes of this work, it is sufficient to consider the very simple Rusanov (local Lax-Friedrichs) numerical flux
\begin{equation} \label{eq:lxf}
	\mathcal{F}(\mathbf{v},\mathbf{w}) = \frac12 \left( \mathbf{f}(\mathbf{v}) + \mathbf{f}(\mathbf{w}) - \alpha (\mathbf{w} - \mathbf{v}) \right).
\end{equation}
Here, $\alpha$ is the parameter of numerical viscosity. In explicit schemes one has to choose $\alpha = \max\{\|\mathbf{f}'(\mathbf{v})\|,\|\mathbf{f}'(\mathbf{w})\|\}$, where $\|\mathbf{f}'(\cdot)\|$ denotes the spectral radius of the Jacobian of the flux function $\mathbf{f}$. In implicit schemes this choice would deserve more attention, especially if the fluid speed is much lower than the sound speed and if one is not interested in the acoustic waves. In this case one can use more appropriate speed estimates in the approximate Riemann solver as explained, e.g., in~\cite{2011DegondTang}. We will present our choice for each test problem. In any case, when assembling the Jacobian of $\mathbf{G}$ of \eqref{eq:predictor:system} for the Newton step \eqref{eq:predicor:newton},
we neglect the terms $\partial_{\mathbf{v}}\alpha$ and $\partial_{\mathbf{w}}\alpha$ in the derivative of $\mathcal{F}$ and approximate it as
\begin{equation*}
	\partial_{\mathbf{v}}\mathcal{F}(\mathbf{v},\mathbf{w}) 
	\approx \tfrac12 \mathrm{J}_{\mathbf{f}}(\mathbf{v}) 
	+ \tfrac12\alpha \mathbb{I}_m 
	\quad\text{and}\quad    
	\partial_{\mathbf{w}}\mathcal{F}(\mathbf{v},\mathbf{w}) 
	\approx \tfrac12 \mathrm{J}_{\mathbf{f}}(\mathbf{w}) 
	- \tfrac12\alpha \mathbb{I}_m,
\end{equation*}
where $\mathrm{J}_{\mathbf{f}}$ denotes the Jacobian of the exact flux function and $\mathbb{I}_m$ is the $m\times m$ identity matrix. For simplicity and not to distract from the main focus of the paper, which is on the implicit time integration of high order schemes, we have used the simple Rusanov flux.

For the time integration, we employ the three stage third order DIRK scheme of~\cite{Alexander1977} with Butcher tableau
\begin{equation*}
	\setlength\arraycolsep{10pt}
	\begin{array}{c|ccc}
		\lambda & \lambda & 0 & 0 \\[1.5ex]
		\frac{(1+\lambda)}{2} & \frac{(1-\lambda)}{2} & \lambda & 0 \\[1.5ex]
		1 & -\frac32 \lambda^2 + 4\lambda - \frac14 & \frac32\lambda^2-5\lambda+\frac54 & \lambda \\[1ex]
		\hline
		&&&\\[-1.5ex]
		& -\frac32 \lambda^2 + 4\lambda - \frac14 & \frac32\lambda^2-5\lambda+\frac54 & \lambda
	\end{array}
\end{equation*}
where $\lambda=0.4358665215$. This method is A-stable and stiffly accurate, and therefore it is L-stable. As a consequence, the Butcher tableau~\eqref{eq:tableau:be} of the composite backward Euler becomes
\begin{equation*}
	\setlength\arraycolsep{10pt}
	\begin{array}{c|ccc}
		\lambda & \lambda & 0 & 0 \\[1.5ex]
		\frac{(1+\lambda)}{2} & \lambda & \frac{1-\lambda}{2} & 0 \\[1.5ex]
		1 & \lambda & \frac{1-\lambda}{2} & \frac{1-\lambda}{2} \\[1ex]
		\hline
		&&&\\[-1.5ex]
		& \lambda & \frac{1-\lambda}{2} & \frac{1-\lambda}{2}
	\end{array}
\end{equation*}
Although this is not necessary for our construction, note that the sequence of the abscissae of this Runge-Kutta schemes is strictly increasing. Other choices of DIRK schemes are possible, e.g.~see~\cite{2009KetchesonMacdonaldGottlieb}.

\begin{table}[t]
	\caption{List of the schemes tested, labels and line type (for those that appear in the plots)}
	\begin{center}
		\begin{tabular}{c|c|c}
			%\hline
			\underline{Label} & \underline{Scheme description} & \underline{Line type} \\[2ex]
			%\hline\hline
			$\Q$ & Quinpi solution without time-limiting & \includegraphics[width=0.1\textwidth]{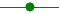} %\begin{tikzpicture}\draw[green!50!black] (1,0)--(2,0);\node[green!50!black] at (1.5,0) {\pgfuseplotmark{*}};\end{tikzpicture} 
			\\[2ex]
			$\QI{1}$ & Quinpi solution with $\mathcal{I}_1$ time-limiting strategy~\eqref{eq:strategy1} & \includegraphics[width=0.1\textwidth]{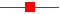}
			%\begin{tikzpicture}\draw[blue] (1,0)--(2,0);\node[blue] at (1.5,0) {\pgfuseplotmark{square*}};\end{tikzpicture}
			\\[2ex]
			$\QI{2}$ & Quinpi solution with $\mathcal{I}_2$ time-limiting strategy~\eqref{eq:strategy2} &
			\\[2ex]
			$\QI{3}$ & Quinpi solution with $\mathcal{I}_3$ time-limiting strategy~\eqref{eq:strategy3} & \includegraphics[width=0.1\textwidth]{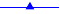}
			%\begin{tikzpicture}\draw[red] (1,0)--(2,0);\node[red] at (1.5,0) {\pgfuseplotmark{triangle*}};\end{tikzpicture}
			\\[2ex]
			%$\Qexp$ & $\Q$ with CFL set as in explicit schemes & \includegraphics[width=0.1\textwidth]{lineType/Q3exp}
			%\begin{tikzpicture}\draw[dashed,green!50!black] (1,0)--(2,0);\node[green!50!black] at (1.5,0) {\pgfuseplotmark{o}};\end{tikzpicture}
			%\\[2ex]
			%$\QIexp{1}$ & $\QI{1}$ with CFL set as in explicit schemes & \includegraphics[width=0.1\textwidth]{lineType/Q3I1exp}
			%\begin{tikzpicture}\draw[dashed,blue] (1,0)--(2,0);\node[blue] at (1.5,0) {\pgfuseplotmark{square}};\end{tikzpicture}
			%\\[2ex]
			%$\QIexp{2}$ & $\QI{2}$ with CFL set as in explicit schemes & 
			%\\[2ex]
			%$\QIexp{3}$ & $\QI{3}$ with CFL set as in explicit schemes & \includegraphics[width=0.1\textwidth]{lineType/Q3I3exp}
			%\begin{tikzpicture}\draw[dashed,red] (1,0)--(2,0);\node[red] at (1.5,0) {\pgfuseplotmark{triangle}};\end{tikzpicture}
			%\\[2ex]
			$\CW$ & Explicit CWENOZ scheme of~\cite{STP23:cweno:boundary} &
			\includegraphics[width=0.1\textwidth]{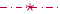} %\begin{tikzpicture}\draw[orange] (1,0)--(2,0);\node[orange] at (1.5,0) {\pgfuseplotmark{star}};\end{tikzpicture}
			\\[2ex]
			$\QP$ & Implicit scheme of~\cite{PSV23:Quinpi} with cell-centered time-limiter & \includegraphics[width=0.1\textwidth]{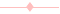} %\begin{tikzpicture}\draw[pink] (1,0)--(2,0);\node[pink] at (1.5,0) {\pgfuseplotmark{diamond*}};\end{tikzpicture}
			\\[2ex]
			$\QPMOOD$ & Implicit scheme of~\cite{VTSP23:Quinpi:Book} which uses MOOD as time-limiter & \includegraphics[width=0.1\textwidth]{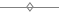} %\begin{tikzpicture}\draw[gray] (1,0)--(2,0);\node[gray] at (1.5,0) {\pgfuseplotmark{diamond}};\end{tikzpicture}
			\\
			%\hline
		\end{tabular}
	\end{center}
	\label{tab:labels}
\end{table}

All numerical simulations are performed with the schemes described in Table~\ref{tab:labels}. For the time-limited schemes, we choose the threshold parameters in \eqref{eq:strategy1} and \eqref{eq:strategy3} as $\gamma_1 = h$ and $\sigma = 10^{-10}$, respectively. Instead, $\gamma_2$, cf.~\eqref{eq:strategy2}, will be presented for each numerical test and we will investigate the performance of the scheme based on its value. In general, the choices of the thresholds $\gamma_1$ and $\gamma_2$ are dictated by the decay of the numerical entropy indicators, cf.~Section~\ref{sec:quinpi:entropy:indicators}, and by the regions of the solution one aims to mark as irregular. Indeed, $\gamma_1=h$ allows us to mark as irregular all regions of the solution where $S^3 > h$, i.e.~those regions characterized by shocks, contact waves, and rarefaction corners. Similarly, a constant value of $\gamma_2$ allows us to distinguish smooth regions, where the indicator $\nicefrac{S^3}{S^1} \sim \mathcal{O}(h^2)$, from irregular regions where $\nicefrac{S^3}{S^1} \sim \mathcal{O}(1)$.

For the solution of the nonlinear systems in the predictor and in the corrector, we used the Newton-Raphson schemes \eqref{eq:predicor:newton} and \eqref{eq:corrector:newton}, and a simple GMRES linear solver with a ILU0 preconditioner. For a detailed discussion on properties of iterative solvers in the context of implicit discretizations of hyperbolic equations we refer to~\cite{2022BirkenLinders}.

\subsection{Gas-dynamics problems} \label{sec:numerics:euler}

In this section, we test the third order Quinpi scheme on several problems based on the one-dimensional nonlinear Euler system for gas-dynamics
\begin{equation} \label{eq:euler:system}
	\partial_t 
	\left( \begin{array}{c}
		\rho \\ \rho v \\ E
	\end{array}\right) +
	\partial_x 
	\left( \begin{array}{c}
		\rho v \\ \rho v^2 + p \\ v(E+p)
	\end{array}\right)  = 0,
\end{equation}
where $\rho$, $v$, $p$ and $E$ are the density, velocity, pressure and energy per unit volume of an ideal gas, whose equation of state is $ E = \frac{p}{\gamma-1} + \frac12 \rho v^2, $ where $\gamma = 1.4$. The eigenvalues of the Jacobian of the flux are $\lambda_1 = v-c$, $\lambda_2=v$ and $\lambda_3=v+c$, where $c=\sqrt{\nicefrac{\gamma p}{\rho}}$ is the sound speed. For this system, we consider the entropy pair defined by the entropy function $\eta(\mathbf{u}) = -\rho \log\left(\frac{p}{(\gamma-1)\rho^\gamma}\right)$ and the entropy flux $\psi(\mathbf{u}) = -v\eta(\mathbf{u})$, see~\cite{GodlewskiRaviart96}. Since we aim to investigate the accuracy of Quinpi schemes on material waves, the parameter of numerical viscosity in~\eqref{eq:lxf} is chosen as $\alpha = \max\{|v^{-}|,|v^{+}|\}$, where $v^{-}$ and $v^{+}$ denote the BED of the velocity of the gas.

\subsubsection{Convergence test}

We check the order of accuracy on the nonlinear system~\eqref{eq:euler:system} by simulating the advection of a density perturbation. The initial condition is
\begin{equation} \label{eq:conv:test:ic}
	\left(\rho,v,p\right) = \left(1+0.5\sin(2\pi x), 1, 10^\kappa\right),
\end{equation}
so that the exact solution at a given time $t$ is a traveling wave
\begin{equation} \label{eq:conv:test:exact}
	\left(\rho,v,p\right) = \left(1+0.5\sin\left(2\pi (x-t)\right), 1, 10^\kappa\right),
\end{equation}
with $\kappa\in\mathbb{N}_0$.
We run the problem up to the final time $t=1$. We use periodic boundary conditions on the domain $[0,1]$ and, thus, we test not only the scheme accuracy, but also the proposed treatment of reconstructions for cells at the boundary.

For this problem, the spectral radius is $\lambda_{\max} := \max_{i=1,2,3} |\lambda_i| = 1 + 10^{\nicefrac{\kappa}{2}}\sqrt{2\gamma}$, where $\lambda_i$ denotes the $i$-th eigenvalue of the system, which is preserved during the time evolution. For an explicit scheme the CFL stability condition would require to choose a grid ratio
$$
\frac{\DT}{h} \leq \frac{1}{\lambda_{\max}} \sim 10^{-\nicefrac{\kappa}{2}}.
$$ 
In order to test the accuracy of the implicit scheme, we run the simulation with the following grid ratio
$$
\frac{\DT}{h} = 4,
$$
and consider $\kappa = 0$ and $\kappa = 4$, corresponding to the Courant numbers $C = \lambda_{\max}\tfrac{\DT}{h} \approx 10.7$ and $C \approx 673.3$, respectively.

The $L^1$- and $L^\infty$-norm errors computed on the density component and the corresponding experimental orders of convergence are listed in Table~\ref{tab:rates:kappa0} and Table~\ref{tab:rates:kappa4} for the two different Courant numbers. The results are obtained with $\gamma_2=0.1$.

While the schemes $\QI{1}$ and $\QI{3}$ achieve the theoretical order of convergence, the scheme $\QI{2}$ exhibits, instead, first order accuracy. This behavior is explained by investigating the time-limiting procedure for $\QI{2}$. In fact, the marking strategy $\mathcal{I}_2$ fails on flat regions where both $S^3$ and $S^1$ are floating point zeros. This causes the unnecessary limiting of all numerical fluxes, and thus the computed solution ends up coinciding with the first order predictor. Instead, the other detectors for time-limiting are able to signal the regularity of the solution, avoiding the unnecessary activation of the time-limiting. For these reasons, in the following we will not consider the scheme $\QI{2}$ anymore.

Furthermore, we observe that Table~\ref{tab:rates:kappa0} and Table~\ref{tab:rates:kappa4} present identical errors and rates of convergence for the methods $\QI{1}$ and $\QI{3}$. This happens because the time-limiting is indeed off on the smooth profile: the two strategies detect no irregular cells and coincide with the unlimited scheme $\Q$. Note that the magnitude of the errors is the same for $\kappa=0$ and $\kappa=4$, even though the latter uses a time-step which violets the explicit CFL by a factor of $60$ with respect to the first one.

%\addtolength{\tabcolsep}{3pt}
\begin{table}[th!]
	\caption{Comparisons of the orders of convergence of Quinpi schemes on the initial condition~\eqref{eq:conv:test:ic} with $\kappa=0$. The grid ratio is $\nicefrac{\Delta t}{h}=4$ which corresponds to a Courant number $C=10.7$.\label{tab:rates:kappa0}}
	\centering
	\vspace{0.25cm}
	\pgfplotstabletypeset[
	font=\small,
	col sep=comma,
	sci zerofill,
	empty cells with={--},
	every head row/.style={before row={\toprule
			&\multicolumn{2}{c}{$\QI{1}$}
			&\multicolumn{2}{c}{$\QI{2}$}
			&\multicolumn{2}{c}{$\QI{3}$}
			\\
		},
		after row=\midrule
	},
	every last row/.style={after row=\bottomrule},
	create on use/rate1/.style={create col/dyadic refinement rate={1}},
	create on use/rate2/.style={create col/dyadic refinement rate={2}},
	create on use/rate3/.style={create col/dyadic refinement rate={3}},
	columns/0/.style={column name={$N$}},
	columns/1/.style={column name={$L^1$ error},sci e},
	columns/rate1/.style={fixed zerofill,column name={rate}},
	columns/2/.style={column name={$L^1$ error},sci e},
	columns/rate2/.style={fixed zerofill,column name={rate}},
	columns/3/.style={column name={$L^1$ error},sci e},
	columns/rate3/.style={fixed zerofill,column name={rate}},
	columns={0,1,rate1,2,rate2,3,rate3},
	skip rows between index={7}{8}
	]
	{press1e0_err1.err}
	\\
	\vspace{0.25cm}
	\pgfplotstabletypeset[
	font=\small,
	col sep=comma,
	sci zerofill,
	empty cells with={--},
	every head row/.style={before row={\toprule
			&\multicolumn{2}{c}{$\QI{1}$}
			&\multicolumn{2}{c}{$\QI{2}$}
			&\multicolumn{2}{c}{$\QI{3}$}
			\\
		},
		after row=\midrule
	},
	every last row/.style={after row=\bottomrule},
	create on use/rate1/.style={create col/dyadic refinement rate={1}},
	create on use/rate2/.style={create col/dyadic refinement rate={2}},
	create on use/rate3/.style={create col/dyadic refinement rate={3}},
	columns/0/.style={column name={$N$}},
	columns/1/.style={column name={$L^\infty$ error},sci e},
	columns/rate1/.style={fixed zerofill,column name={rate}},
	columns/2/.style={column name={$L^\infty$ error},sci e},
	columns/rate2/.style={fixed zerofill,column name={rate}},
	columns/3/.style={column name={$L^\infty$ error},sci e},
	columns/rate3/.style={fixed zerofill,column name={rate}},
	columns={0,1,rate1,2,rate2,3,rate3},
	skip rows between index={7}{8}
	]
	{press1e0_errinf.err}
\end{table}

\begin{table}[th!]
	\caption{Comparisons of the orders of convergence of Quinpi schemes on the initial condition~\eqref{eq:conv:test:ic} with $\kappa=4$. The grid ratio is $\nicefrac{\Delta t}{h}=4$ which corresponds to a Courant number $C=673.3$.\label{tab:rates:kappa4}}
	\centering
	\vspace{0.25cm}
	\pgfplotstabletypeset[
	font=\small,
	col sep=comma,
	sci zerofill,
	empty cells with={--},
	every head row/.style={before row={\toprule
			&\multicolumn{2}{c}{$\QI{1}$}
			&\multicolumn{2}{c}{$\QI{2}$}
			&\multicolumn{2}{c}{$\QI{3}$}
			\\
		},
		after row=\midrule
	},
	every last row/.style={after row=\bottomrule},
	create on use/rate1/.style={create col/dyadic refinement rate={1}},
	create on use/rate2/.style={create col/dyadic refinement rate={2}},
	create on use/rate3/.style={create col/dyadic refinement rate={3}},
	columns/0/.style={column name={$N$}},
	columns/1/.style={column name={$L^1$ error},sci e},
	columns/rate1/.style={fixed zerofill,column name={rate}},
	columns/2/.style={column name={$L^1$ error},sci e},
	columns/rate2/.style={fixed zerofill,column name={rate}},
	columns/3/.style={column name={$L^1$ error},sci e},
	columns/rate3/.style={fixed zerofill,column name={rate}},
	columns={0,1,rate1,2,rate2,3,rate3}
	]
	{press1e4_err1.err}
	\\
	\vspace{0.25cm}
	\pgfplotstabletypeset[
	font=\small,
	col sep=comma,
	sci zerofill,
	empty cells with={--},
	every head row/.style={before row={\toprule
			&\multicolumn{2}{c}{$\QI{1}$}
			&\multicolumn{2}{c}{$\QI{2}$}
			&\multicolumn{2}{c}{$\QI{3}$}
			\\
		},
		after row=\midrule
	},
	every last row/.style={after row=\bottomrule},
	create on use/rate1/.style={create col/dyadic refinement rate={1}},
	create on use/rate2/.style={create col/dyadic refinement rate={2}},
	create on use/rate3/.style={create col/dyadic refinement rate={3}},
	columns/0/.style={column name={$N$}},
	columns/1/.style={column name={$L^\infty$ error},sci e},
	columns/rate1/.style={fixed zerofill,column name={rate}},
	columns/2/.style={column name={$L^\infty$ error},sci e},
	columns/rate2/.style={fixed zerofill,column name={rate}},
	columns/3/.style={column name={$L^\infty$ error},sci e},
	columns/rate3/.style={fixed zerofill,column name={rate}},
	columns={0,1,rate1,2,rate2,3,rate3}
	]
	{press1e4_errinf.err}
\end{table}

\subsubsection{Stiff Riemann problems} \label{sec:numerics:stiffRP}

In the following, we test Quinpi schemes on stiff Riemann problems for the Euler system~\eqref{eq:euler:system}, namely for those problems characterized by $\frac{|v|}{|v|+c} \ll 1, \quad \forall\,(x,t)$.
In this case, we expect a strong separation between acoustic and material waves, because they travel with speeds having different magnitudes. In particular, the acoustic speeds would impose a strict constraint on the time-step
$$
\frac{\DT_{\text{stab}}}{h} \leq \frac{1}{\max_{x} |v|+c}.
$$
With implicit schemes, instead, the time-step is not dictated by stability and therefore we will choose it in order to reproduce accurately the slow material wave:
$$
\frac{\DT_{\text{acc}}}{h} = \frac{1}{|v_{\text{cw}}|},
$$
where $v_{\text{cw}}$ denotes an estimate of the velocity of material waves. The ratio between the time-step $\DT_{\text{stab}}$ required for stability and the time-step $\DT_{\text{acc}}$ required for accuracy on the material wave, is a measure of the stiffness of the problem. We define the Courant number as
$$
C = \frac{\DT_{\text{acc}}}{\Delta t_{\text{stab}}}.
$$

\begin{figure}
	\centering
	\begin{subfigure}[b]{0.32\textwidth}
		\centering
		\includegraphics[width=\textwidth]{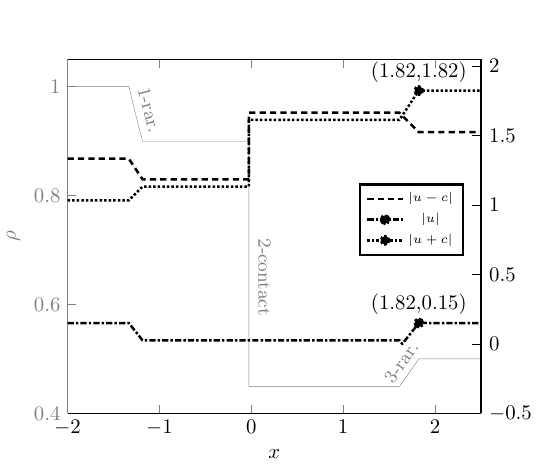}
		\caption{Non-symmetric expansion.}
		\label{fig:exact_A}
	\end{subfigure}
	\hfill
	\begin{subfigure}[b]{0.32\textwidth}
		\centering
		\includegraphics[width=\textwidth]{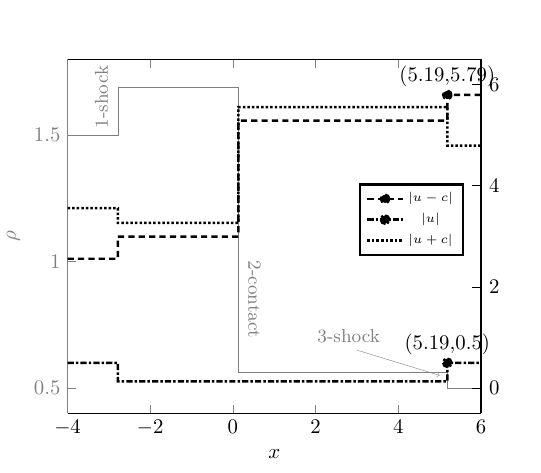}
		\caption{Colliding flows.}
		\label{fig:exact_B}
	\end{subfigure}
	\hfill
	\begin{subfigure}[b]{0.32\textwidth}
		\centering
		\includegraphics[width=\textwidth]{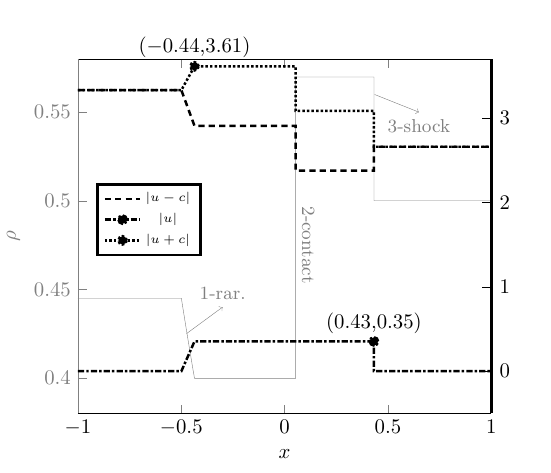}
		\caption{Stiff Lax.}
		\label{fig:exact_C}
	\end{subfigure}
	\caption{Exact density solution (thin gray solid line) and wave speeds (thick lines) of the Riemann problems at final time.}
	\label{fig:exact}
\end{figure}

We consider the following Riemann problems.

\begin{description}
	\item[Test (a): Non-symmetric expansion problem.]
	\begin{equation} \label{eq:tworar:ic}
		\left(\rho,u,p\right) = \begin{cases}
			(1, -0.15, 1), & x < 0\\
			(0.5,0.15,1), & x \geq 0.
		\end{cases}
	\end{equation}
	The exact density at time $t=1$ is depicted in the left panel of Figure~\ref{fig:exact} and is characterized by a 1-rarefaction moving with negative velocity, a 3-rarefaction moving with positive velocity, and a 2-contact wave with velocity $v_{\text{cw}}=-2.57\times10^{-2}$.
	For this problem, one has
	$
	\max_{(x,t)} |v|+c = 1.82, %\ \text{ and } \ 3.8\cdot10^{-4} \leq \mathrm{Ma} \leq 0.17.
	$
	thus an explicit scheme would require to choose a time-step $\DT_{\text{stab}}$ such that
	$$
	\frac{\DT_{\text{stab}}}{h} \leq 0.549.
	$$
	With an implicit scheme, it is possible to overcome this stability restriction. In particular, choosing a time-step $\DT_{\text{acc}}$ to meet the accuracy on the contact wave, one has
	$$
	\frac{\DT_{\text{acc}}}{h} = 6.66,
	$$
	which implies a Courant number $C \approx 12.1$.
	
	\item[Test (b): colliding flow problem.]
	\begin{equation} \label{eq:twoshocks:ic}
		\left(\rho,v,p\right) = \begin{cases}
			(1.5,  0.5, 10), & x < 0\\
			(0.5, -0.5, 10), & x \geq 0.
		\end{cases}
	\end{equation}
	The exact density at time $t=1$ is depicted in the center panel of Figure~\ref{fig:exact} and is characterized by a 1-shock wave moving with negative velocity, a 3-shock wave moving with positive velocity, and a 2-contact wave with velocity $v_{\text{cw}}=0.13$.
	For this problem, one has
	$
	\max_{(x,t)} |v|+c = 5.79, %\ \text{ and } \ 2.5\cdot10^{-2} \leq \mathrm{Ma} \leq 0.16.
	$
	thus an explicit scheme would require to choose a time-step $\DT_{\text{stab}}$ such that
	$$
	\frac{\DT_{\text{stab}}}{h} \leq 0.17.
	$$
	With an implicit scheme, instead, one can choose
	$$
	\frac{\DT_{\text{acc}}}{h} = 2,
	$$
	which implies a Courant number $C \approx 11.6$.
	
	\item[Test (c): modified Lax shock tube problem.] \begin{equation} \label{eq:stiffLax:ic}
		\left(\rho,v,p\right) = \begin{cases}
			(0.445,  0, 3.528), & x < 0\\
			(0.5, 0, 2.528), & x \geq 0.
		\end{cases}
	\end{equation}
	Compared to the classical Lax shock tube problem, the initial condition is characterized by a zero left velocity and a higher right pressure. This allows a faster separation of the waves. The exact density at time $t=0.15$ is depicted in the right panel of Figure~\ref{fig:exact} and is characterized by a 1-rarefaction moving with negative velocity, a 3-shock moving with positive velocity, and a 2-contact wave with velocity $v_{\text{cw}}=0.35$.
	For this problem, one has
	$
	\max_{(x,t)} |v|+c = 3.61,
	$
	thus an explicit scheme would require to choose a time-step $\DT_{\text{stab}}$ such that
	$$
	\frac{\DT_{\text{stab}}}{h} \leq 0.28.
	$$
	With an implicit scheme, instead, one can choose
	$$
	\frac{\DT_{\text{acc}}}{h} = 2.83,
	$$
	which implies a Courant number $C \approx 10.2$.
\end{description}

\begin{figure}[t!]
	\centering
	\begin{subfigure}[b]{0.4\textwidth}
		\centering
		\includegraphics[width=0.7\textwidth]{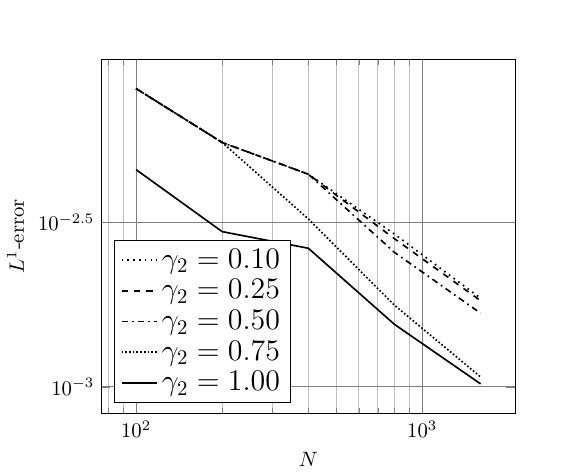}
		\caption{Non-symmetric expansion at $t=1$.}
		\label{fig:gamma2_A}
	\end{subfigure}
	\hfill
	\begin{subfigure}[b]{0.28\textwidth}
		\centering
		\includegraphics[width=\textwidth]{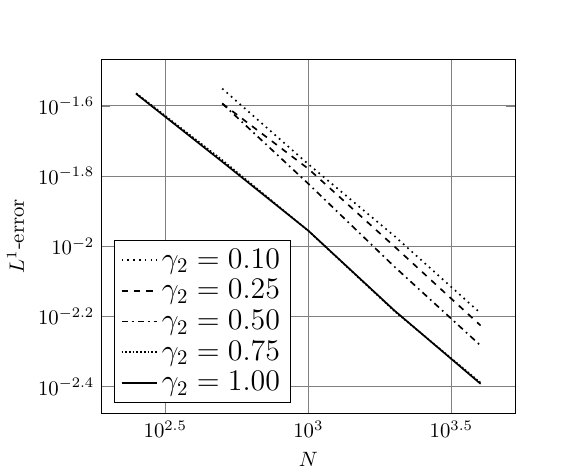}
		\caption{Colliding flows at $t=1$.}
		\label{fig:gamma2_B}
	\end{subfigure}
	\hfill
	\begin{subfigure}[b]{0.28\textwidth}
		\centering
		\includegraphics[width=\textwidth]{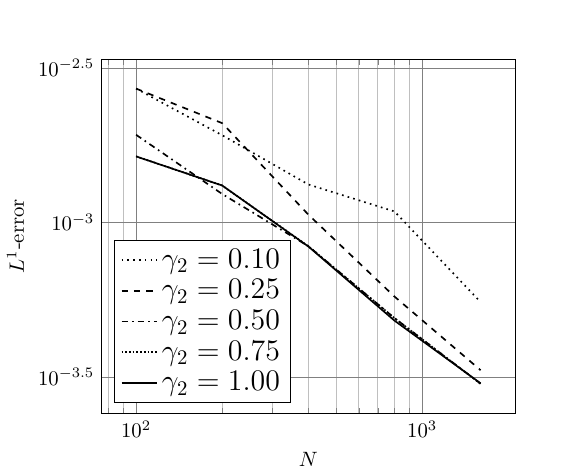}
		\caption{Stiff Lax at $t=0.15$.}
		\label{fig:gamma2_C}
	\end{subfigure}
	\caption{$L^1$-errors on the contact wave computed with the $\QI{3}$ schemes on several grids and different values of the threshold parameter $\gamma_2$.}
	\label{fig:gamma2}
\end{figure}

Here below we comment the results of all three tests. First, we investigate the role of the parameter threshold $\gamma_2$, see~\eqref{eq:strategy3}, which defines the time-limiting strategy $\mathcal{I}_3$. In particular, we study the accuracy of the resulting $\QI{3}$ scheme on the contact wave for several values of $\gamma_2$. In Figure~\ref{fig:gamma2}, we show the $L^1$-errors as function of the number of cells. We observe that in all the proposed Riemann problems the choice $\gamma_2=1$ provides smaller errors. Therefore, all the numerical solutions computed with $\QI{3}$ are considered with $\gamma_2=1$.

\begin{figure}[t!]
	\centering
	\begin{subfigure}[b]{\textwidth}
		\centering
		\includegraphics[width=\textwidth]{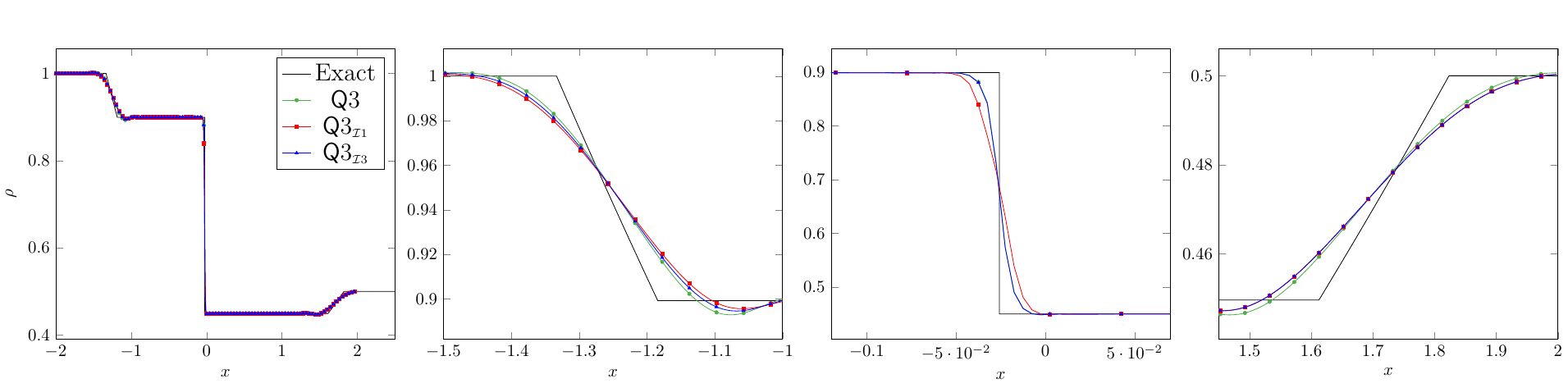}
		\caption{Non-symmetric expansion at $t=1$ with $N=800$ cells.}
		\label{fig:timelim_A}
	\end{subfigure}
	%\hfill
	\begin{subfigure}[b]{\textwidth}
		\centering
		\includegraphics[width=\textwidth]{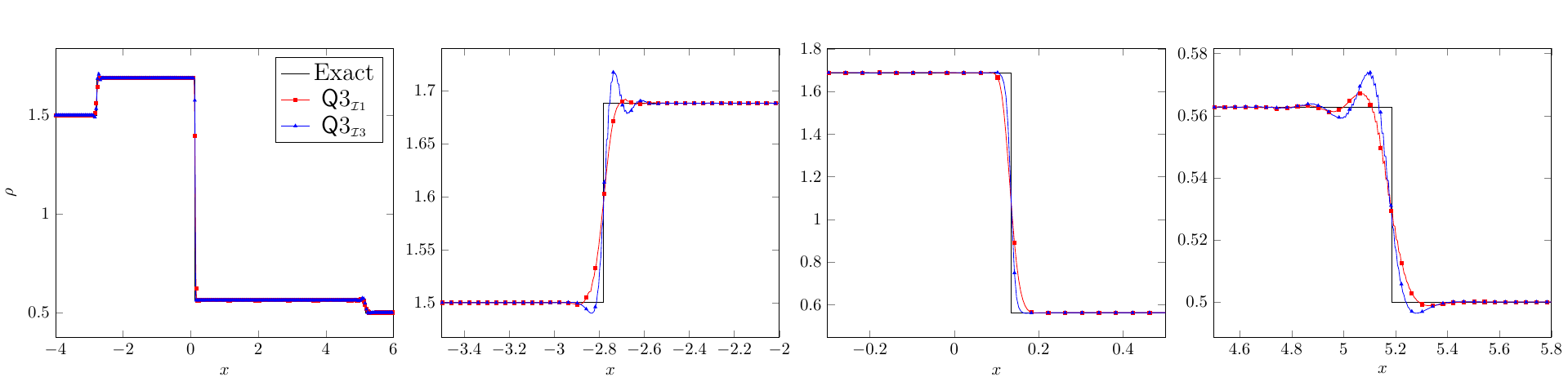}
		\caption{Colliding flows at $t=1$ with $N=2000$ cells.}
		\label{fig:timelim_B}
	\end{subfigure}
	%\hfill
	\begin{subfigure}[b]{\textwidth}
		\centering
		\includegraphics[width=\textwidth]{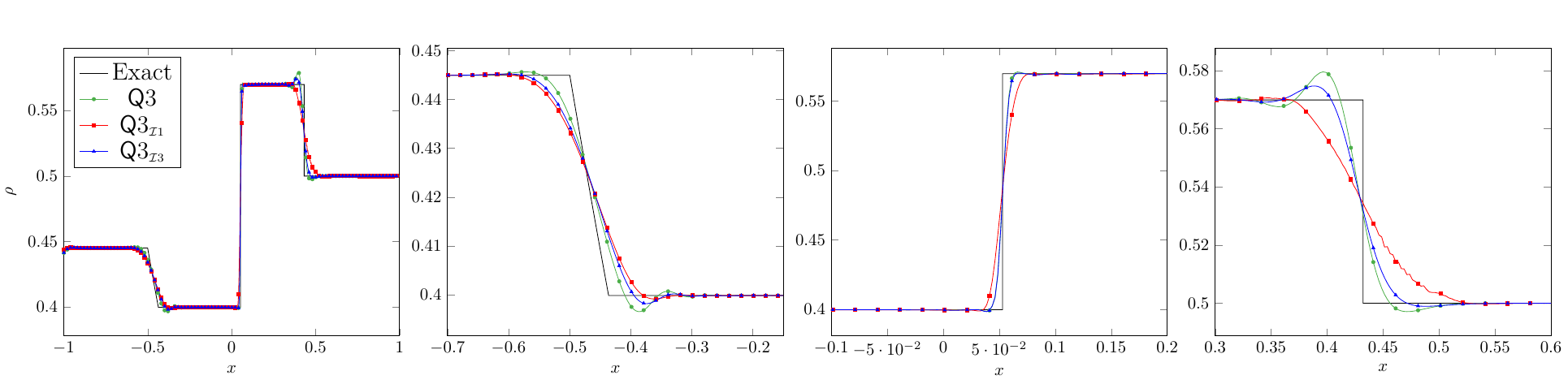}
		\caption{Stiff Lax at $t=0.15$ with $N=800$ cells.}
		\label{fig:timelim_C}
	\end{subfigure}
	\caption{Comparison between third order implicit schemes, without time-limiting and with different time-limiting strategies. The left-most panel shows the density solution. The other panels show zooms on the three waves. The $\Q$ solution of the colliding flows problem is not showed because it blows up before reaching the final time.}
	\label{fig:timelim}
\end{figure}

\begin{figure}[t!]
	\centering
	\begin{subfigure}[b]{0.32\textwidth}
		\centering
		\includegraphics[width=\textwidth]{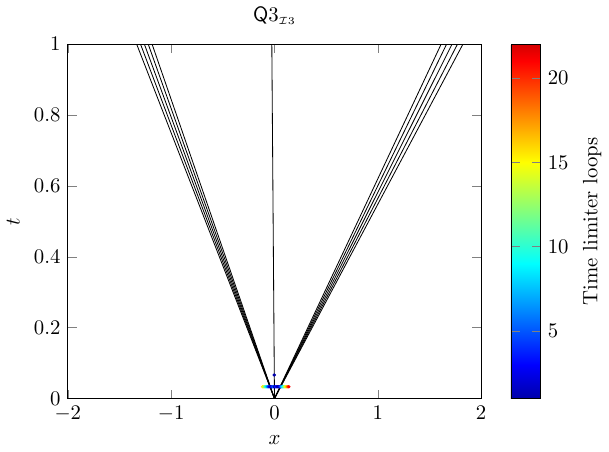}
		\caption{Non-symmetric expansion\\with $N=800$ cells.}
		\label{fig:xt_timelim_A}
	\end{subfigure}
	%\hfill
	\begin{subfigure}[b]{0.32\textwidth}
		\centering
		\includegraphics[width=\textwidth]{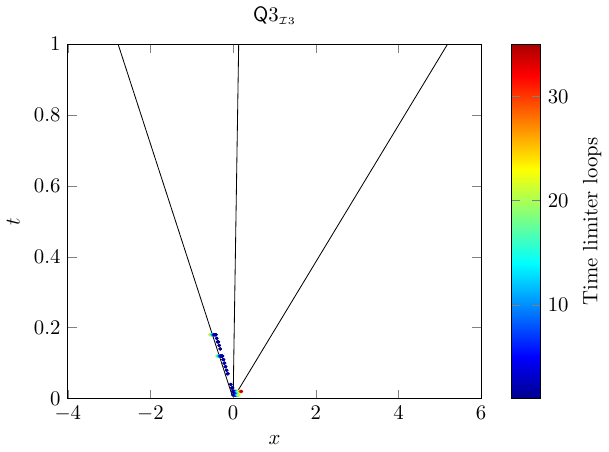}
		\caption{Colliding flows with\\$N=2000$ cells.}
		\label{fig:xt_timelim_B}
	\end{subfigure}
	%\hfill
	\begin{subfigure}[b]{0.32\textwidth}
		\centering
		\includegraphics[width=\textwidth]{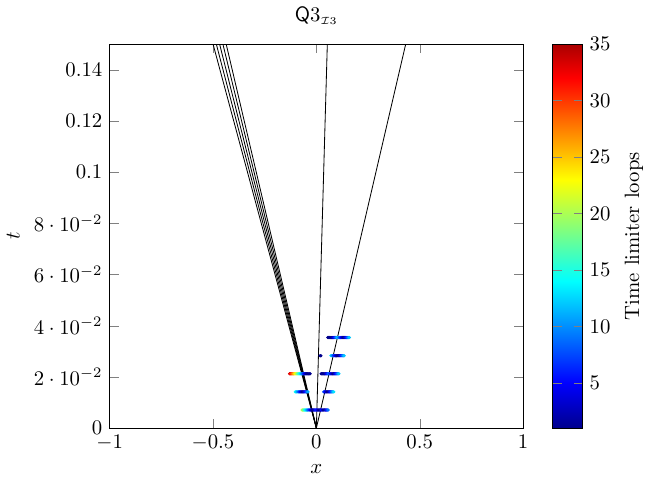}
		\caption{Stiff Lax with\\$N=800$ cells.}
		\label{fig:xt_timelim_C}
	\end{subfigure}
	\caption{Structure of the solutions in the $x$-$t$ diagram. Limited cell interfaces are marked at the time levels and locations where time-limiting is applied. The colors denote the loop number at which a cell interface is marked: interfaces limited at initial loops are marked in blue, while those limited at final loops are marked in red. These results can be compared with the data provided in Table~\ref{tab:stat}.}
	\label{fig:xt:timelim}
\end{figure}

\begin{table}[t!]
	\caption{Statistics on the time-limiting procedure $\mathcal{I}_3$ of Quinpi schemes. These data are graphically depicted in Figure~\ref{fig:xt:timelim}.\label{tab:stat}}
	\begin{subtable}[h]{0.32\textwidth}
		\caption{Non-symmetric expansion.\label{tab:stat:tworar}}
		\centering
		%\vspace{0.25cm}
		\pgfplotstabletypeset[
		font=\scriptsize,
		col sep=comma,
		sci zerofill,
		empty cells with={--},
		every head row/.style={before row={\toprule
				&\multicolumn{2}{c}{$\QI{3}$}%&\multicolumn{2}{c}{$\QIexp{3}$}
				\\[2ex]
			},
			after row=\midrule
		},
		every last row/.style={after row=\bottomrule},
		columns/0/.style={column name={$N$}},
		columns/4/.style={column name={\shortstack{max.~number of \\ limited fluxes}}},
		columns/5/.style={column name={\shortstack{\% limited \\ time-steps}}},
		%columns/9/.style={column name={\shortstack{max.~number of \\ limited fluxes}}},
		%columns/10/.style={column name={\shortstack{\% limited \\ time-steps}}},
		columns={0,4,5}%,9,10},
	%skip rows between index={0}{2}
	]
	{tworar_ind313_stat.out}
\end{subtable}
\hfill
\begin{subtable}[h]{0.3175\textwidth}
	%	\vspace{0.25cm}
	\caption{Colliding flows.\label{tab:stat:twoshocks}}
	\centering
	%\vspace{0.25cm}
	\pgfplotstabletypeset[
	font=\scriptsize,
	col sep=comma,
	sci zerofill,
	empty cells with={--},
	every head row/.style={before row={\toprule
			&\multicolumn{2}{c}{$\QI{3}$}%&\multicolumn{2}{c}{$\QIexp{3}$}
			\\[2ex]
		},
		after row=\midrule
	},
	every last row/.style={after row=\bottomrule},
	columns/0/.style={column name={$N$}},
	columns/4/.style={column name={\shortstack{max.~number of \\ limited fluxes}}},
	columns/5/.style={column name={\shortstack{\% limited \\ time-steps}}},
	%columns/9/.style={column name={\shortstack{max.~number of \\ limited fluxes}}},
	%columns/10/.style={column name={\shortstack{\% limited \\ time-steps}}},
	columns={0,4,5}%,9,10},
%skip rows between index={0}{2}
]
{twoshocks_ind313_stat.out}
\end{subtable}
\hfill
\begin{subtable}[h]{0.3175\textwidth}
%\vspace{0.25cm}
\caption{Stiff Lax.\label{tab:stat:stiffLax}}
\centering
%\vspace{0.25cm}
\pgfplotstabletypeset[
font=\scriptsize,
col sep=comma,
sci zerofill,
empty cells with={--},
every head row/.style={before row={\toprule
		&\multicolumn{2}{c}{$\QI{3}$}%&\multicolumn{2}{c}{$\QIexp{3}$}
		\\[2ex]
	},
	after row=\midrule
},
every last row/.style={after row=\bottomrule},
columns/0/.style={column name={$N$}},
columns/4/.style={column name={\shortstack{max.~number of \\ limited fluxes}}},
columns/5/.style={column name={\shortstack{\% limited \\ time-steps}}},
%columns/9/.style={column name={\shortstack{max.~number of \\ limited fluxes}}},
%columns/10/.style={column name={\shortstack{\% limited \\ time-steps}}},
columns={0,4,5}%,9,10},
%skip rows between index={0}{2}
]
{stifflax_ind313_stat.out}
\end{subtable}
\end{table}

Next, we show the need for time-limiting. In Figure~\ref{fig:timelim} we compare the solutions computed without time-limiting ($\Q$) and with the two limiting strategies $\mathcal{I}_1$ and $\mathcal{I}_3$. The figures clearly show that
\begin{itemize}
\item time-limiting is needed to avoid spurious oscillations. In particular, the numerical solution of the colliding flows problem obtained with the $\Q$ scheme is not showed because it blows up at time $t=0.93$;
\item strategy $\mathcal{I}_1$, see~\eqref{eq:strategy1}, is more diffusive. Therefore, in the following we will consider the third strategy $\mathcal{I}_3$ only, see~\eqref{eq:strategy3}. Moreover, observe that the $\QI{3}$ solution is accurate as the $\Q$ solution on the contact waves of the non-symmetric expansion and of the stiff Lax problems, while reducing the oscillations on the acoustic waves;
\item the time-limiting is active on all the waves at initial times, i.e.~when they separate themselves, cf.~Figure~\ref{fig:xt:timelim}. At larger times, instead, limiting acts in particular on shock waves but stops when the waves are sufficiently separated. In Table~\ref{tab:stat} we observe that the percentage of limited steps is larger in problems where shock waves develop. Moreover, the number of limited fluxes is large because during the large time-step a wave propagates across more cells in a single time-step.
\end{itemize}

\begin{figure}[t!]
\centering
\begin{subfigure}[b]{\textwidth}
\centering
\includegraphics[width=\textwidth]{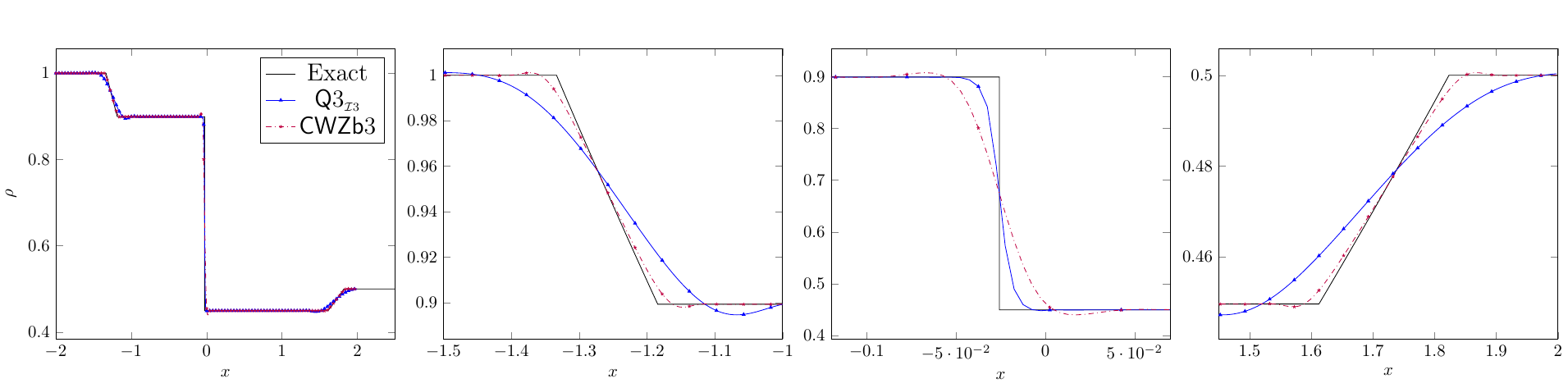}
\caption{Non-symmetric expansion at $t=1$ with $N=800$ cells. The explicit solution is obtained with time-step $2.7\times10^{-3}$, whereas the implicit solution is computed with time-step $3.3\times10^{-2}$.}
\label{fig:imex_A}
\end{subfigure}
%\hfill
\begin{subfigure}[b]{\textwidth}
\centering
\includegraphics[width=\textwidth]{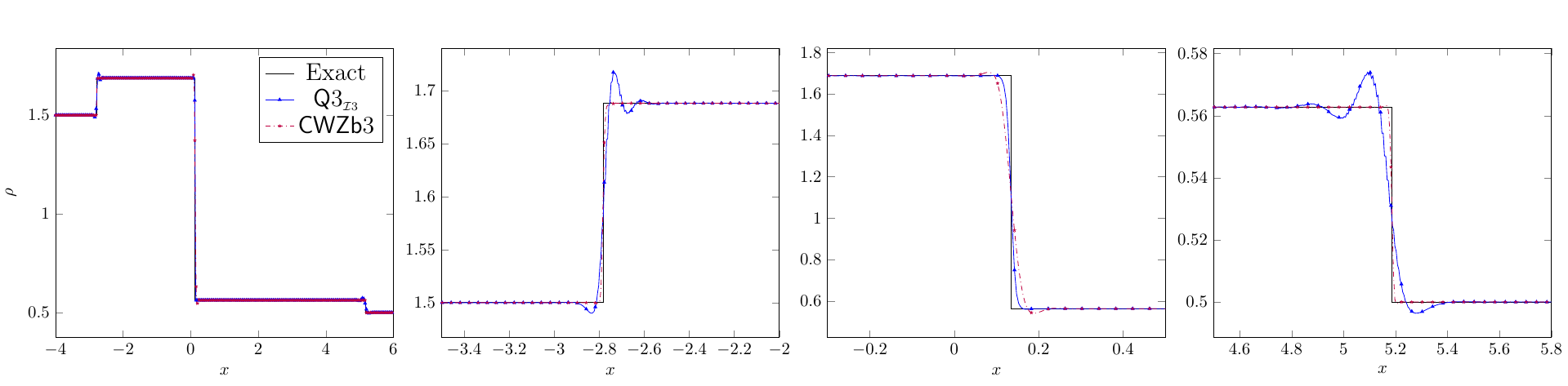}
\caption{Colliding flows at $t=1$ with $N=2000$ cells. The explicit solution is obtained with time-step $8.5\times10^{-4}$, whereas the implicit solution is computed with time-step $1.0\times10^{-2}$.}
\label{fig:imex_B}
\end{subfigure}
%\hfill
\begin{subfigure}[b]{\textwidth}
\centering
\includegraphics[width=\textwidth]{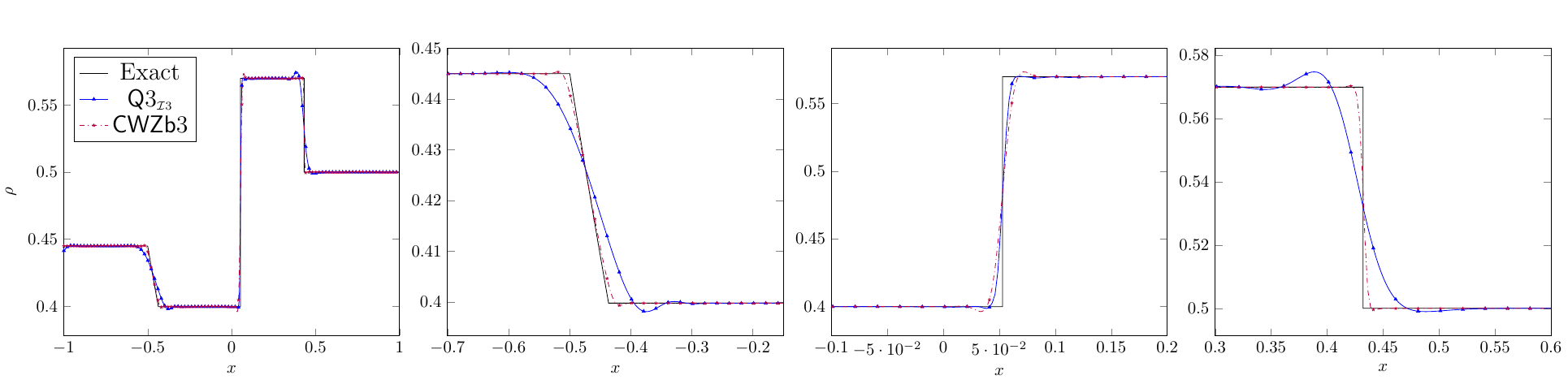}
\caption{Stiff Lax at $t=0.15$ with $N=800$ cells. The explicit solution is obtained with time-step $7.0\times10^{-4}$, whereas the implicit solution is computed with time-step $7.1\times10^{-3}$.}
\label{fig:imex_C}
\end{subfigure}
\caption{Comparison between third order implicit and explicit schemes. The left-most panel shows the density solution. The other panels show zooms  on the three waves.}
\label{fig:imex}
\end{figure}

Finally, we compare the explicit scheme $\CW$ of~\cite{STP23:cweno:boundary}, which has to be run with $\DT_{\text{stab}}$, and the implicit scheme $\QI{3}$ run at $\DT_{\text{acc}}$, i.e.~centered on the material wave.
Figure~\ref{fig:imex} shows that
\begin{itemize}
\item as expected, the explicit scheme is more accurate than the implicit one on the acoustic waves;
\item on the contact wave, which is the one we are focusing on, the implicit scheme is more accurate than the explicit one and runs with a time-step that is about $10$ times larger.
\end{itemize}

\subsection{Low-Mach problems} \label{sec:numerics:lowMach}

In this section, we test the performance of Quinpi schemes on low-Mach problems, namely such that
$$
\mathrm{Ma} := \frac{|v|}{c} \ll 1, \quad \forall\,(x,t).
$$
Here, the aim is to show that the implicit framework we propose has the ability to guarantee stability and accuracy on such problems. We do not claim that the Quinpi approach is asymptotic preserving, e.g.~see~\cite{2017DimarcoLoubereVignal}.

In order to describe the low-Mach number limit, we consider the rescaled compressible Euler system
\begin{equation} \label{eq:euler:system:rescaled}
\partial_t 
\left( \begin{array}{c}
\rho \\ \rho v \\ E
\end{array}\right) +
\partial_x 
\left( \begin{array}{c}
\rho v \\ \rho v^2 + \frac{1}{\varepsilon^2} p \\ v(E+p)
\end{array}\right)  = 0,
\end{equation}
with the equation of state $E = \frac{p}{\gamma-1} + \frac{\varepsilon^2}{2} \rho v^2$, $\gamma=1.4$.
The parameter $\varepsilon$ describes the Mach number within the non-dimensionalized system, and one has $\mathrm{Ma} = \nicefrac{\varepsilon}{\sqrt{\gamma}}$. System~\eqref{eq:euler:system:rescaled} is hyperbolic with eigenvalues $\lambda_1 = v - \nicefrac{c}{\varepsilon}$, $\lambda_2 = v$, $\lambda_3= v + \nicefrac{c}{\varepsilon}$.
For a more detailed discussion on the low-Mach number scaling we refer, e.g., to~\cite{2018BoscarinoRussoScandurra,1995Klein}. Although in the following we consider only smooth problems, we still choose the parameter of numerical viscosity in~\eqref{eq:lxf} according to the material speed, namely $\alpha = \max\{|v^{-}|,|v^{+}|\}$, where $v^{-}$ and $v^{+}$ denote the BED of the velocity of the gas.

\subsubsection{Convergence test}

We compute the experimental order of convergence of the Quinpi scheme on the computational domain [-2.5,2.5] using the smooth initial condition~\cite{2018BoscarinoRussoScandurra}
$$
\rho(x,0) = \left( 1 + \varepsilon \frac{(\gamma-1) u(x,0)}{2\sqrt{\gamma}} \right)^{\frac{2}{\gamma-1}}, \quad u(x,0) = \sin\left(\frac{2\pi x}{5}\right), \quad p(x,0) = \rho(x,0)^\gamma.
$$
The order of accuracy of the scheme is investigated at final time $t=0.3$ for $\varepsilon_1 = 0.8$, $\varepsilon_2=0.3$, and at final time $t=0.01$ for $\varepsilon_3=10^{-4}$. For this problem, $|v| \leq 1$, $\forall\,(x,t)$, and $\lambda_{\max}(\varepsilon) := \max_{(x,t)} |v|+\nicefrac{c}{\varepsilon}$ is
$$
\lambda_{\max}(0.8) = 2.6786, \quad \lambda_{\max}(0.3)=5.1439, \quad \lambda_{\max}(10^{-4})=1.1833 \times 10^{4}.
$$
For an explicit scheme the CFL stability condition would require to choose a grid ratio
\begin{align*}
\frac{\DT_{\varepsilon_1}}{h} &\leq \frac{1}{\lambda_{\max}(0.8)} \approx 0.3733, \\ \frac{\DT_{\varepsilon_2}}{h} &\leq \frac{1}{\lambda_{\max}(0.3)} \approx 0.1944, \\ \frac{\DT_{\varepsilon_3}}{h} &\leq \frac{1}{\lambda_{\max}(10^{-4})} \approx 8.45\times10^{-5}.
\end{align*}
We overcome the CFL condition running the simulation using the following grid ratios:
$$
\frac{\DT_{\varepsilon_i}}{h} = \frac{C}{\lambda_{\max}(\varepsilon_i)}, \quad i=1,2,3,
$$
with Courant number $C = 20$.

\begin{table}[th!]
\caption{Experimental order of convergence of Quinpi scheme on the low-Mach problem with Courant number $C=20$.\label{tab:ratesLM:cfl20}}
\centering
\vspace{0.25cm}
\pgfplotstabletypeset[
font=\small,
col sep=comma,
sci zerofill,
empty cells with={--},
every head row/.style={before row={\toprule
&\multicolumn{2}{c}{$\varepsilon=0.8$}
&\multicolumn{2}{c}{$\varepsilon=0.3$}
&\multicolumn{2}{c}{$\varepsilon=10^{-4}$}
\\
},
after row=\midrule
},
every last row/.style={after row=\bottomrule},
create on use/rate1/.style={create col/dyadic refinement rate={1}},
create on use/rate2/.style={create col/dyadic refinement rate={2}},
create on use/rate3/.style={create col/dyadic refinement rate={3}},
columns/0/.style={column name={$N$}},
columns/1/.style={column name={$L^1$ error},sci e},
columns/rate1/.style={fixed zerofill,column name={rate}},
columns/2/.style={column name={$L^1$ error},sci e},
columns/rate2/.style={fixed zerofill,column name={rate}},
columns/3/.style={column name={$L^1$ error},sci e},
columns/rate3/.style={fixed zerofill,column name={rate}},
columns={0,1,rate1,2,rate2,3,rate3},
skip rows between index={0}{3}
]
{err1.err}
\\
\vspace{0.25cm}
\pgfplotstabletypeset[
font=\small,
col sep=comma,
sci zerofill,
empty cells with={--},
every head row/.style={before row={\toprule
&\multicolumn{2}{c}{$\varepsilon=0.8$}
&\multicolumn{2}{c}{$\varepsilon=0.3$}
&\multicolumn{2}{c}{$\varepsilon=10^{-4}$}
\\
},
after row=\midrule
},
every last row/.style={after row=\bottomrule},
create on use/rate1/.style={create col/dyadic refinement rate={1}},
create on use/rate2/.style={create col/dyadic refinement rate={2}},
create on use/rate3/.style={create col/dyadic refinement rate={3}},
columns/0/.style={column name={$N$}},
columns/1/.style={column name={$L^\infty$ error},sci e},
columns/rate1/.style={fixed zerofill,column name={rate}},
columns/2/.style={column name={$L^\infty$ error},sci e},
columns/rate2/.style={fixed zerofill,column name={rate}},
columns/3/.style={column name={$L^\infty$ error},sci e},
columns/rate3/.style={fixed zerofill,column name={rate}},
columns={0,1,rate1,2,rate2,3,rate3},
skip rows between index={0}{3}
]
{errinf.err}
\end{table}

Density errors for the different values of the Mach number are listed in Table~\ref{tab:ratesLM:cfl20} and computed with the Quinpi scheme $\Q$. The theoretical third order of convergence is observed.

\subsubsection{Two colliding acoustic pulses}

We consider the two colliding acoustic pulses taken from~\cite{1995Klein}. The initial condition is given by
\begin{align*}
\rho(x,0) &= \rho_0 + \frac{\varepsilon\rho_1}{2} \left( 1-\cos\left( \frac{2\pi x}{L} \right) \right), \quad \rho_0=0.955, \quad \rho_1=2, \\
u(x,0) &= -\frac{u_0}{2} \mathrm{sign}(x) \left( 1-\cos\left( \frac{2\pi x}{L} \right) \right), \quad u_0=2\sqrt{\gamma}, \\
p(x,0) &= p_0 + \frac{\varepsilon p_1}{2} \left( 1-\cos\left( \frac{2\pi x}{L} \right) \right), \quad p_0=1, \quad p_1=2\gamma,
\end{align*}
on the computational domain $[-L,L]$, $L=\nicefrac{2}{\varepsilon}$. We use periodic boundary conditions and, thus, we also test the proposed treatment of~\cite{STP23:cweno:boundary} which does not require ghost cells.

We choose $\varepsilon_1=\nicefrac{1}{11}$ and $\varepsilon_2=10^{-4}$. In both cases, $\max_{(x,t)} |v| = 2.3656$, whereas $\lambda_{\max}(\varepsilon_1) = 16.038$ and $\lambda_{\max}(\varepsilon_2) = 1.21\times10^4$. Therefore, using an explicit scheme would require to choose a time-step
$$
\frac{\DT_{\varepsilon_1}}{h} \leq 6.24\times 10^{-2}, \quad \frac{\DT_{\varepsilon_2}}{h} \leq 8.26 \times 10^{-5}.
$$
With an implicit scheme, instead, we can choose
$$
\frac{\DT_{\varepsilon_1}}{h} = \frac{1}{\max_{(x,t)} |v|} = 4.23 \times 10^{-1},
$$
which corresponds to a Courant number $C = 6.78$. We employ the same Courant number for $\varepsilon_2$ which determines
$$
\frac{\DT_{\varepsilon_2}}{h} = 5.59 \times 10^{-4}.
$$

\begin{figure}[t!]
\centering
\begin{subfigure}[b]{\textwidth}
\centering
\includegraphics[width=\textwidth]{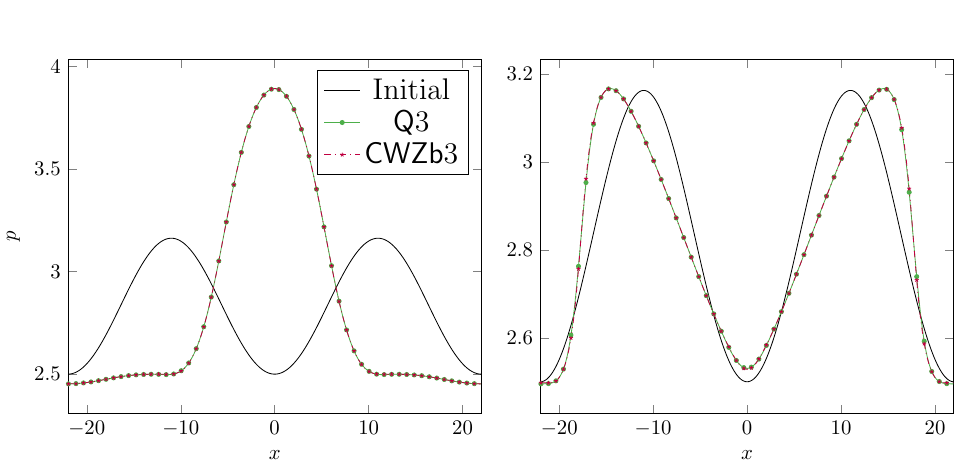}
\caption{$\varepsilon=\frac{1}{11}$, $t=0.815$ (left) and $t=1.63$ (right). The explicit solution is obtained with time-step $0.624$, whereas the implicit solution is computed with time-step $4.23$.}
\label{fig:colliding:pulses:eps1}
\end{subfigure}
%\hfill
\begin{subfigure}[b]{\textwidth}
\centering
\includegraphics[width=\textwidth]{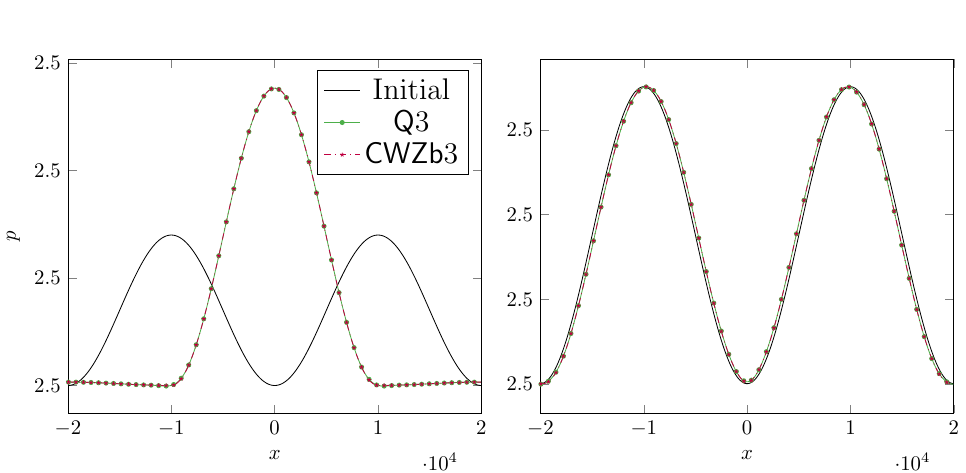}
\caption{$\varepsilon=10^{-4}$, $t=0.815$ (left) and $t=1.63$ (right). The explicit solution is obtained with time-step $7.5 \times 10^{-4}$, whereas the implicit solution is computed with time-step $5.1 \times 10^{-3}$.}
\label{fig:colliding:pulses:eps2}
\end{subfigure}
\caption{Two colliding pulses problem with $n=440$ cells. The solutions show the initial pressure and the final pressure profiles obtained with the explicit $\CW$ and the implicit $\Q$ third order schemes.}
\label{fig:colliding:pulses}
\end{figure}

In Figure~\ref{fig:colliding:pulses} we show the pressure profile obtained using the third order Quinpi scheme $\Q$ without time-limiting at final times $t = 0.815$ and $t = 1.63$, with $N=440$ cells. The solution is compared with the third order explicit scheme $\CW$.

\subsection{Scalar problems} \label{sec:numerics:scalar}

Finally, we provide some numerical simulations on scalar conservation laws. In fact, although this work is mostly focused on hyperbolic systems, the Quinpi framework introduced in~\cite{PSV23:Quinpi}, and further extended in~\cite{VTSP23:Quinpi:Book}, was based on very different time-limiting procedures. Therefore, this section aims to investigate the performance of the proposed MOOD limiter combined with the numerical entropy production indicator on standard scalar problems. In all the scalar problems, the numerical viscosity in~\eqref{eq:lxf} is chosen as $\alpha = \max\{|f'(v)|,|f'(w)|\}$ where $f$ is the flux function. All simulations are performed with $\gamma_2 = 0.1$.

\subsubsection{Linear transport}

We consider the linear scalar conservation law
\begin{equation} \label{eq:linadv}
\partial_t u(x,t) + \partial_x u(x,t) = 0,
\end{equation}
for $(x,t)\in[-1,1]\times(0,2]$, with periodic boundary conditions in space and discontinuous initial conditions:
\begin{subequations}\label{eq:nonsmoothIC}
\begin{equation} \label{eq:sindiscontIC}
u_0(x) = \sin(\pi x) +
\begin{cases}
3, & -0.4 \leq x \leq 0.4,\\
0, & \text{otherwise,}	
\end{cases}
\end{equation}
\begin{equation} \label{eq:doublestepIC}
u_0(x) =
\begin{cases}
1, & -0.25 \leq x \leq 0.25,\\
0, & \text{otherwise.}
\end{cases}
\end{equation}
\end{subequations}
Here, we consider the entropy function $\eta(u) = \frac{u^2}{2}$ and the entropy flux $\psi(u) = \frac{u^2}{2}$.

\begin{figure}[t!]
\centering
\begin{subfigure}[b]{\textwidth}
\centering
\includegraphics[width=\textwidth]{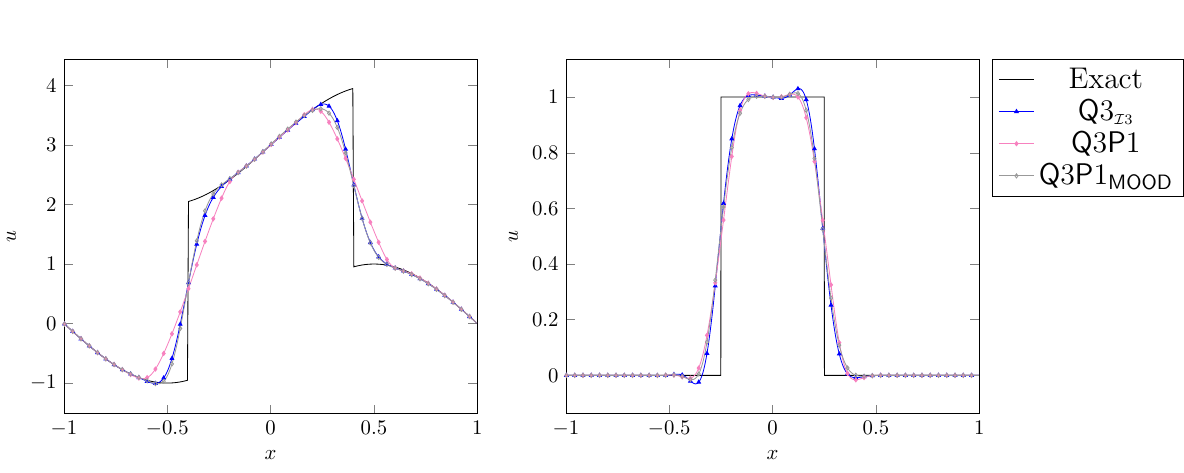}
\caption{$\frac{\Delta t}{h}=5$.}
\label{fig:lintra:cou5}
\end{subfigure}
%\hfill
\begin{subfigure}[b]{\textwidth}
\centering
\includegraphics[width=\textwidth]{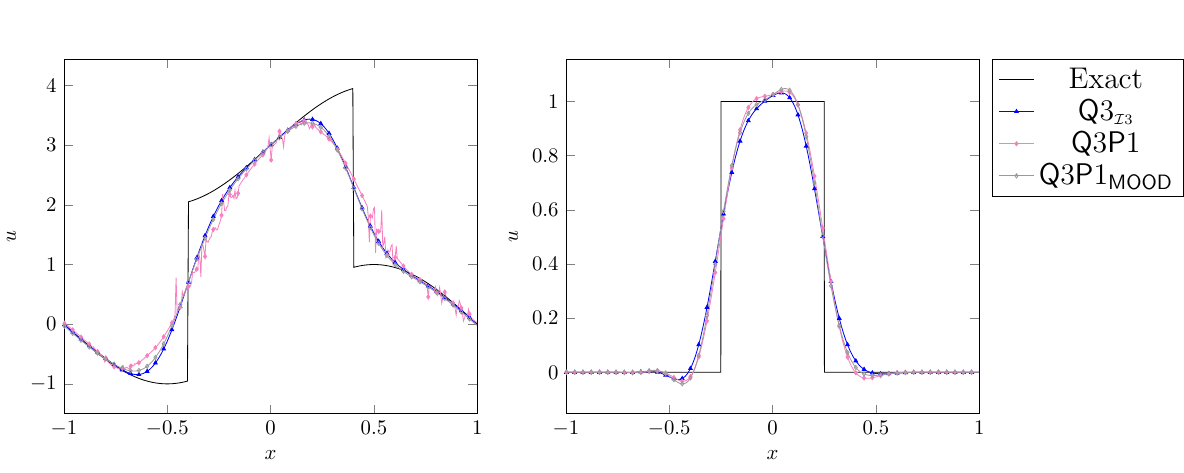}
\caption{$\frac{\Delta t}{h}=10$.}
\label{fig:lintra:cou10}
\end{subfigure}
\caption{Linear transport equation~\eqref{eq:linadv} with initial condition~\eqref{eq:sindiscontIC} (left panels) and initial condition~\eqref{eq:doublestepIC} (right panels). All the solution are obtained on a grid of $400$ cells with two Courant numbers. The label $\QP$ refers to the scheme of~\cite{PSV23:Quinpi} with mass redistribution.}
\label{fig:lintra}
\end{figure}

The results are reported in Figure~\ref{fig:lintra}, computed with two different Courant numbers, i.e.~$\nicefrac{\Delta t}{h}=5$ and $\nicefrac{\Delta t}{h}=10$ on $400$ space cells. As introduced in Table~\ref{tab:labels}, the solution labeled by $\QP$ refers to the scheme of~\cite{PSV23:Quinpi}, whereas the solution labeled by $\QPMOOD$ refers to the one of~\cite{VTSP23:Quinpi:Book}. We note that the scheme $\QI{3}$ performs better than the other two implicit schemes on the discontinuous sinusoidal profile since it produces lower dissipation around the jump discontinuities. However, the low-dissipation properties leads to lightly larger oscillations around the discontinuities of the double-step profile. In particular, with $\nicefrac{\Delta t}{h} = 10$, the upper flat part of the solution is not reproduced by all the schemes.

\subsubsection{Burgers' equation}

Next, we compare the schemes on the nonlinear Burgers' equation
\begin{equation} \label{eq:burgers}
\partial_t u(x,t) + \partial_x \left( \frac{u^2(x,t)}{2} \right) = 0,
\end{equation}
for $x\in[-1,1]$, with periodic boundary conditions in space. Here, we consider the entropy function $\eta(u) = \frac{u^2}{2}$ and the entropy flux $\psi(u) = \frac{u^3}{3}$.

%\begin{figure}[t!]
%\centering
%\includegraphics[width=\textwidth]{burgers_doublestep}
%\caption{Burgers' equation~\eqref{eq:burgers} with initial condition~\eqref{eq:doublestepIC} on $400$ cells at time $t=0.5$, with $\nicefrac{\Delta t}{h} = 3$ (left), $\nicefrac{\Delta t}{h} = 5$ (middle), $\nicefrac{\Delta t}{h} = 10$ (right).}
%\label{fig:burgers:doublestep}
%\end{figure}

%First, we consider the discontinuous initial condition~\eqref{eq:doublestepIC} and compute the numerical solution at time $t=0.5$ on $400$ space cells with three different Courant numbers, i.e.~$\nicefrac{\Delta t}{h}=3,5,10$. The exact solution is characterized by a standing-tail with right-moving head rarefaction and a shock with positive velocity. The results are reported in Figure~\ref{fig:burgers:doublestep}. We observe that the scheme $\QPMOOD$ reproduces the rarefaction and the shock with less dissipation when $\nicefrac{\Delta t}{h}=3$. However, the time-limiting of $\QPMOOD$ leads to higher dissipation at increasing Courant numbers. Moreover, while the dissipation of $\QP$ is lower than the dissipation of the other schemes on the rarefaction head at $\nicefrac{\Delta t}{h}=3$, we notice that $\QI{3}$ reproduces accurately the tail of the rarefaction and avoids the undershoot on the shock wave.

%Furthermore,
We test%the Burgers' equation
~\eqref{eq:burgers} on the smooth initial condition
\begin{equation} \label{eq:shockinteractionIC}
u_0(x) = 0.2 -\sin(\pi x) + \sin(2 \pi x), 
\end{equation}
computing the solution at three different times, namely $t=\nicefrac{1}{2\pi},0.6,1$, and Courant numbers $\nicefrac{\Delta t}{h}=3,10$, on a space grid of $400$ cells. The exact solution is characterized by the formation of two shocks which collide developing a single discontinuity.

\begin{figure}[t!]
\centering
\begin{subfigure}[b]{\textwidth}
\centering
\includegraphics[width=\textwidth]{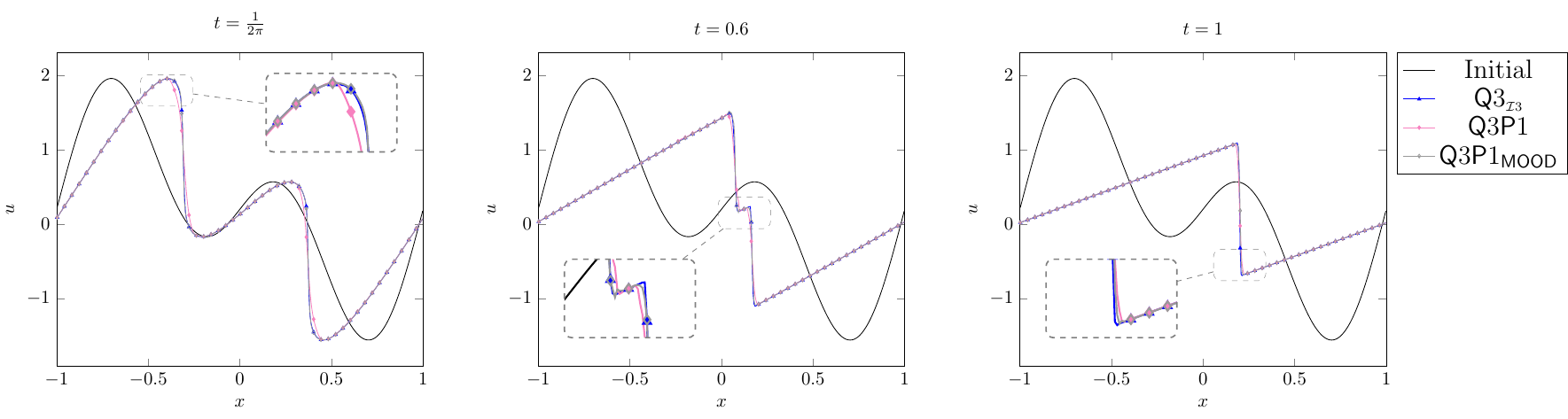}
\caption{$\frac{\Delta t}{h}=3$.}
\label{fig:burgers:shockinteraction:cou3}
\end{subfigure}
%\hfill
\begin{subfigure}[b]{\textwidth}
\centering
\includegraphics[width=\textwidth]{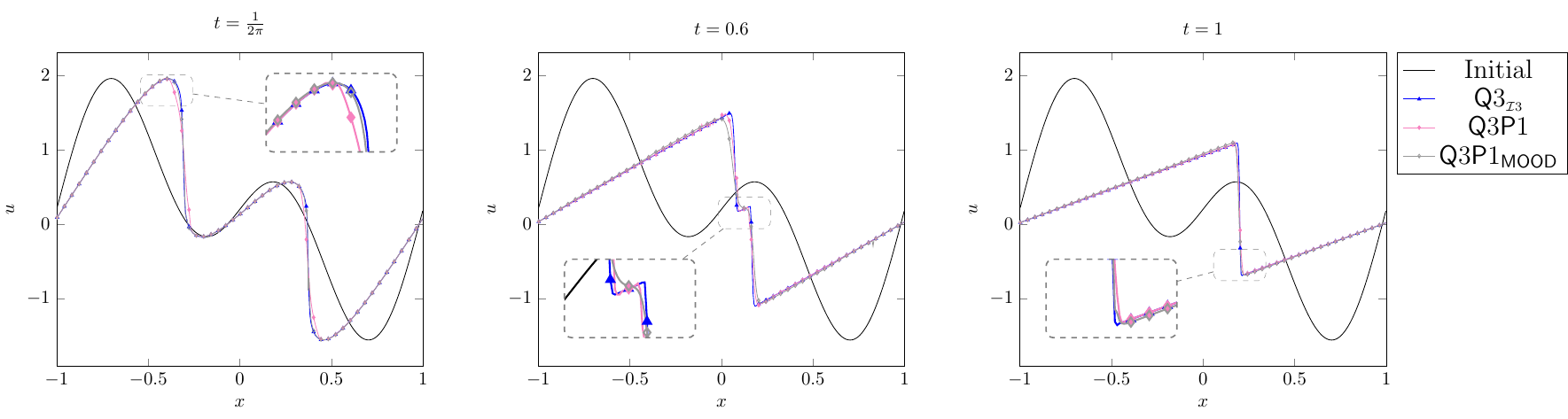}
\caption{$\frac{\Delta t}{h}=10$.}
\label{fig:burgers:shockinteraction:cou10}
\end{subfigure}
\caption{Burgers' equation~\eqref{eq:burgers} with initial condition~\eqref{eq:shockinteractionIC} on $400$ cells with $\nicefrac{\Delta t}{h}=3$ (top panels) and $\nicefrac{\Delta t}{h}=10$ (bottom panels), at three different times.}
\label{fig:burgers:shockinteraction}
\end{figure}

The numerical results are depicted in Figure~\ref{fig:burgers:shockinteraction}. We observe that all schemes do not produce spurious oscillations at time $t=1/(2\pi)$, just before the two shocks appear, but the scheme $\QP$ is more diffusive. At $t=0.6$, $\QPMOOD$ exhibits a small undershoot on the left shock, clearly visible in the zoom box, which is not present in the solution computed by $\QI{3}$. Increasing the Courant number produces more dissipation within the numerical schemes $\QP$ and $\QPMOOD$, whereas $\QI{3}$ is more accurate and does not exhibit spurious oscillations.

\section{Conclusions and future research} \label{sec:conclusion}

In this work, we have proposed a framework for the high order implicit integration of hyperbolic systems of conservation laws. The approach is based on the recent work~\cite{PSV23:Quinpi}, which was focused on scalar conservation laws.
The crucial idea of Quinpi is to use a first order implicit predictor to freeze the nonlinear space limiters in the CWENO third order reconstruction, because such a predictor gives a pre-estimate of the WENO weights, which will be used by the higher order method. Space limiting however is not enough, and it is necessary to perform limiting in time, based on mixing the computed high order solution and the low order predictor. As a consequence, the nonlinearity of the scheme is linked only to the nonlinearity of the flux, while space and time-limiting need not enter into the nonlinear iterations needed to update the solution.

At variance with~\cite{PSV23:Quinpi}, we have designed a flux-centered conservative time-limiting procedure inspired by the MOOD technique~\cite{CDL11:MOOD,CDL12:MOOD}, where the numerical entropy residual~\cite{Puppo04:numerical:entropy,PS11:numerical:entropy} has been used as indicator of non-smooth solutions and troubled cells. This is more robust than the previous cell-centered limiter that required a conservative correction after the time-limiting.

We mention that this implicit approach has not been tailored to the solution of a specific problem. We expect the new scheme to have impact in many applications described by stiff problems, in particular in those cases in which the phenomena of interest travel much slower than the fastest wave speeds. Our tests show that implicit schemes can provide a better accuracy compared to explicit ones on material waves, even running with a time-step about 10 times larger.

Future work will involve the study of the Quinpi approach combined with different time integration methods, for instance using hybrid BDF schemes, extensions to higher order, and possibly investigation of analytical properties of the scheme. 
Further attention will also be paid to improving the efficiency of the scheme. To this end, we will investigate ad-hoc linear solvers and preconditioners, especially in view of higher dimensional applications, and we plan to improve the time-limiting procedure. %A preliminary analysis on the work-precision shows that the accuracy of the implicit and explicit schemes is similar on the material wave. Whereas, the implicit scheme becomes competitive, and more convenient than the explicit method, when the Courant number $C \gtrsim 12$.
We expect to make the computational cost scale more favorably for the implicit approach.

\paragraph{Acknowledgments} This work was partly supported by MUR (Ministry of University and Research) under the PRIN-2017 project on ``Innovative Numerical Methods for Evolutionary Partial Differential Equations and Applications'' (number 2017KKJP4X), and under the PRIN-2022 project on ``High order structure-preserving semi-implicit schemes for hyperbolic equations'' (number 2022JH87B4).

This work was carried out within the Ateneo Sapienza projects 2022 ``Approssimazione numerica di modelli differenziali e applicazioni'' and 2023 ``Modeling, numerical treatment of hyperbolic equations and optimal control problems''.

GP and GV acknowledge support from PNRR-MUR project ``Italian Research Center on High Performance Computing, Big Data and Quantum Computing'' and from PNRR-MUR project PE0000013-FAIR. 

The authors are members of the INdAM Research National Group of Scientific Computing (INdAM-GNCS).

The authors appreciate the unknown referee’s valuable and profound comments which helped to improve significantly the results of this work.

No conflict of interest is extant in the present work.

\bibliographystyle{plain}
\bibliography{Quinpi}

\begin{thebibliography}{10}

\bibitem{2017AbbateAllSpeed}
E.~Abbate, A.~Iollo, and G.~Puppo.
\newblock An all-speed relaxation scheme for gases and compressible materials.
\newblock {\em J. Comput. Phys.}, 351:1--24, 2017.

\bibitem{ABC16:improvedWENOZ}
F.~Acker, R.~B.~de~R. Borges, and B.~Costa.
\newblock An improved {WENO-Z} scheme.
\newblock {\em J. Comput. Phys.}, 313:726--753, 2016.

\bibitem{Alexander1977}
R.~Alexander.
\newblock Diagonally implicit {R}unge-{K}utta methods for stiff {O.D.E.}'s.
\newblock {\em {SIAM} J. Numer. Anal.}, 14(6):1006--1021, 1977.

\bibitem{2020Arbogast}
T.~Arbogast, C.S. Huang, X.~Zhao, and D.~N. King.
\newblock A third order, implicit, finite volume, adaptive {R}unge-{K}utta
  {WENO} scheme for advection-diffusion equations.
\newblock {\em Comput. Methods Appl. Mech. Engrg.}, 368, 2020.

\bibitem{ABIR:19}
S.~Avgerinos, F.~Bernard, A.~Iollo, and G.~Russo.
\newblock Linearly implicit all {M}ach number shock capturing schemes for the
  {E}uler equations.
\newblock {\em J. Comput. Phys.}, 393:278--312, 2019.

\bibitem{Balsara:AOWENO}
D.~S. Balsara, S.~Garain, and C.~W. Shu.
\newblock An efficient class of {WENO} schemes with adaptive order.
\newblock {\em J. Comput. Phys.}, 326:780--804, 2016.

\bibitem{BB:23}
W.~Barsukow and R.~Borsche.
\newblock Implicit active flux methods for linear advection.
\newblock Preprint arXiv:2303.13318, 2023.

\bibitem{Marsha:AMR}
M.~J. Berger and R.~J. Le~Veque.
\newblock Adaptive mesh refinement using wave-propagation algorithms for
  hyperbolic systems.
\newblock {\em {SIAM} J. Numer. Anal.}, 35(6):2298--2316, 1998.

\bibitem{2022BirkenLinders}
P.~Birken and V.~Linders.
\newblock Conservation properties of iterative methods for implicit
  discretizations of conservation laws.
\newblock {\em J. Sci. Comp.}, 92(60), 2022.

\bibitem{BoscarinoQiuRussoXiong}
S.~Boscarino, J.~Qiu, G.~Russo, and T.~Xiong.
\newblock High order semi-implicit {WENO} schemes for all-{M}ach full {E}uler
  system of gas dynamics.
\newblock {\em {SIAM} J. Sci. Comput.}, 44(2):B368--B394, 2022.

\bibitem{2018BoscarinoRussoScandurra}
S.~Boscarino, G.~Russo, and L.~Scandurra.
\newblock All {M}ach number second order semi-implicit scheme for the {E}uler
  equations of gas dynamics.
\newblock {\em J. Sci. Comp.}, 77(2):850--884, 2018.

\bibitem{CCD:11}
M.~Castro, B.~Costa, and W.~S. Don.
\newblock High order weighted essentially non-oscillatory {WENO-Z} schemes for
  hyperbolic conservation laws.
\newblock {\em J. Comput. Phys.}, 230(5):1766--1792, 2011.

\bibitem{CDL11:MOOD}
S.~Clain, S.~Diot, and R.~Loub\`{e}re.
\newblock A high-order finite volume method for hyperbolic systems:
  {M}ulti-dimensional {O}ptimal {O}rder {D}etection ({MOOD}).
\newblock {\em J. Comput. Phys.}, 230(10):4028--4050, 2011.

\bibitem{CDL12:MOOD}
S.~Clain, S.~Diot, and R.~Loub\`{e}re.
\newblock Improved detection criteria for the {M}ulti-dimensional {O}ptimal
  {O}rder {D}etection {MOOD} on unstructured meshes with very high-order
  polynomials.
\newblock {\em Comp. \& Fluids}, 64:43--63, 2012.

\bibitem{CNPT:09}
F.~Coquel, Q.~L. Nguyen, M.~Postel, and Q.~H. Tran.
\newblock Local time stepping with adaptive time step control for a two-phase
  fluid system.
\newblock In {\em ESAIM: {P}roceedings}, volume~29, pages 73--88, 2009.

\bibitem{CNPT:10}
F.~Coquel, Q.~L. Nguyen, M.~Postel, and Q.~H. Tran.
\newblock Entropy-satisfying relaxation method with large time-steps for
  {E}uler {IBVP}s.
\newblock {\em Math. Comp.}, 79:1493--1533, 2010.

\bibitem{CNPT:10a}
F.~Coquel, Q.~L. Nguyen, M.~Postel, and Q.~H. Tran.
\newblock Local time stepping applied to implicit-explicit methods for
  hyperbolic systems.
\newblock {\em Multiscale Model. Simul.}, 8(2):540--570, 2010.

\bibitem{CPPT:06}
F.~Coquel, M.~Postel, N.~Poussineau, and Q.~H. Tran.
\newblock Multiresolution technique and explicit–implicit scheme for
  multicomponent flows.
\newblock {\em J. Numer. Math.}, 14(3):187--216, 2006.

\bibitem{CPSV:cweno}
I.~Cravero, G.~Puppo, M.~Semplice, and G.~Visconti.
\newblock {CWENO}: uniformly accurate reconstructions for balance laws.
\newblock {\em Math. Comp.}, 87(312):1689--1719, 2018.

\bibitem{CSV19:cwenoz}
I.~Cravero, M.~Semplice, and G.~Visconti.
\newblock Optimal definition of the nonlinear weights in multidimensional
  {C}entral {WENOZ} reconstructions.
\newblock {\em {SIAM} J. Numer. Anal.}, 57(5):2328--2358, 2019.

\bibitem{2011DegondTang}
P.~Degond and M.~Tang.
\newblock All speed scheme for the low {M}ach number limit of the isentropic
  {E}uler equations.
\newblock {\em Comm. Computat. Phys.}, 10(1):1--31, 2011.

\bibitem{2010DellacherieLowMach}
S.~Dellacherie.
\newblock Analysis of {G}odunov type schemes applied to the compressible
  {E}uler system at low {M}ach number.
\newblock {\em J. Comput. Phys.}, 229(4):978--1016, 2010.

\bibitem{2017DimarcoLoubereVignal}
G.~Dimarco, R.~Loub\`{e}re, and M.-H. Vignal.
\newblock Study of a new asymptotic preserving scheme for the {E}uler system in
  the low {M}ach number limit.
\newblock {\em {SIAM} J. Sci. Comput.}, 39(5):2099--2128, 2017.

\bibitem{2007DurasaisamyBaeder}
K.~Duraisamy and J.D. D.~Baeder.
\newblock Implicit scheme for hyperbolic conservation laws using non
  oscillatory reconstruction in space and time.
\newblock {\em {SIAM} J. Sci. Comput.}, 29:2607--2620, 2007.

\bibitem{2003DurasaisamyBaeder}
K.~Duraisamy, J.D. D.~Baeder, and J.~G. Liu.
\newblock Concepts and application of time-limiters to high resolution schemes.
\newblock {\em J. Sci. Comp.}, 19:139--162, 2003.

\bibitem{EBS22:Implicit:Networks}
M.~Eimer, R.~Borsche, and N.~Siedow.
\newblock Implicit finite volume method with a posteriori limiting for
  transport networks.
\newblock {\em Adv. Comput. Math.}, 48(3):21, 2022.

\bibitem{FKRZ:22}
P.~Frolkovi{\v{c}}, S.~Kri{\v{s}}kov\'{a}, M.~Rohov\'{a}, and {\v{Z}}erav\'{y}.
\newblock Semi-implicit methods for advection equations with explicit forms of
  numerical solution.
\newblock {\em Jpn. J. Ind. Appl. Math.}, 39:843--867, 2022.

\bibitem{FZ:22}
P.~Frolkovi{\v{c}} and {\v{Z}}erav\'{y}.
\newblock High resolution compact implicit numerical scheme for conservation
  laws.
\newblock {\em Appl. Math. Comput.}, 442:127720, 2023.

\bibitem{2023GersterSemplice}
S.~Gerster and M.~Semplice.
\newblock Semi-conservative high order scheme with numerical entropy indicator
  for intrusive formulations of hyperbolic systems.
\newblock {\em J. Comput. Phys.}, 489:112254, 2023.

\bibitem{GodlewskiRaviart96}
E.~Godlewski and P.~A. Raviart.
\newblock {\em Numerical approximation of hyperbolic systems of conservation
  laws}, volume 118 of {\em Applied Mathemtical Sciences}.
\newblock Springer-Verlag, New York, 1996.

\bibitem{Gottlieb:iWENO:2006}
S.~Gottlieb, J.~S. Mullen, and S.~J. Ruuth.
\newblock {A Fifth Order Flux Implicit WENO Method}.
\newblock {\em J. Sci. Comp.}, 27:271--287, 2006.

\bibitem{2001Gottlieb}
S.~Gottlieb, C.W. Shu, and E.~Tadmor.
\newblock Strong stability preserving high-order time discretization methods.
\newblock {\em SIAM Rev.}, 43:73--85, 2001.

\bibitem{JiangShu:96}
G.-S. Jiang and C.-W. Shu.
\newblock Efficient implementation of weighted {ENO} schemes.
\newblock {\em J. Comput. Phys.}, 126:202--228, 1996.

\bibitem{2009KetchesonMacdonaldGottlieb}
D.~I. Ketcheson, C.~B. MacDonald, and S.~Gottlieb.
\newblock Optimal implicit strong stability preserving {R}unge-{K}utta methods.
\newblock {\em Appl. Numer. Math.}, 59(2):373--392, 2009.

\bibitem{Ketcheson:fluxbased:2013}
D.~I. Ketcheson, C.~B. MacDonald, and S.~J. Ruuth.
\newblock Spatially partitioned embedded {R}unge-{K}utta methods.
\newblock {\em {SIAM} J. Numer. Anal.}, 51(5):2887--2910, 2013.

\bibitem{1995Klein}
R.~Klein.
\newblock Semi-implicit extension of a {G}odunov-type scheme based on low
  {M}ach number asymptotics. {I}: {O}ne-dimensional flow.
\newblock {\em J. Comput. Phys.}, 121(2):213--237, 1995.

\bibitem{LeVeque:book}
R.~Le~Veque.
\newblock {\em Finite Volume Methods for Hyperbolic Problems}.
\newblock Cambridge Texts in Applied Mathematics. Cambridge University Press,
  2004.

\bibitem{LPR:00:SIAMJSciComp}
D.~Levy, G.~Puppo, and G.~Russo.
\newblock Compact central {WENO} schemes for multidimensional conservation
  laws.
\newblock {\em {SIAM} J. Sci. Comput.}, 22(2):656--672, 2000.

\bibitem{LDD:14}
R.~Loub\`{e}re, M.~Dumbser, and S.~Diot.
\newblock A new family of high order unstructured mood and {ADER} finite volume
  schemes for multidimensional systems of hyperbolic conservation laws.
\newblock {\em Comm. Computat. Phys.}, 16:718--763, 2014.

\bibitem{Puppo04:numerical:entropy}
G.~Puppo.
\newblock Numerical entropy production for central schemes.
\newblock {\em {SIAM} J. Sci. Comput.}, 25(4):1382--1415, 2004.

\bibitem{PS11:numerical:entropy}
G.~Puppo and M.~Semplice.
\newblock Numerical entropy and adaptivity for finite volume schemes.
\newblock {\em Comm. Computat. Phys.}, 10(5):1132--1160, 2011.

\bibitem{PS16:entropy:balance}
G.~Puppo and M.~Semplice.
\newblock Well-balanced high order 1{D} schemes on non-uniform grids and
  entropy residuals.
\newblock {\em J. Sci. Comp.}, 66:1052--1076, 2016.

\bibitem{PSV23:Quinpi}
G.~Puppo, M.~Semplice, and G.~Visconti.
\newblock Quinpi: {I}ntegrating conservation laws with {CWENO} implicit
  methods.
\newblock {\em Comm.~Appl.~Math.~\& Comput.}, 5:343--369, 2023.

\bibitem{SCR:CWENOquadtree}
M.~Semplice, A.~Coco, and G.~Russo.
\newblock Adaptive mesh refinement for hyperbolic systems based on third-order
  compact {WENO} reconstruction.
\newblock {\em J. Sci. Comp.}, 66:692--724, 2016.

\bibitem{SL:18:AMRMOOD}
M.~Semplice and R.~Loub\`{e}re.
\newblock Adaptive-{M}esh-{R}efinement for hyperbolic systems of conservation
  laws based on a posteriori stabilized high order polynomial reconstructions.
\newblock {\em J. Comput. Phys.}, 354:86--110, 2018.

\bibitem{STP23:cweno:boundary}
M.~Semplice, E.~Travaglia, and G.~Puppo.
\newblock One- and multi-dimensional {CWENOZ} reconstructions for implementing
  boundary conditions without ghost cells.
\newblock {\em Comm.~Appl.~Math.~\& Comput.}, 5:143--169, 2023.

\bibitem{SempliceVisconti:2020}
M.~Semplice and G.~Visconti.
\newblock Efficient implementation of adaptive order reconstructions.
\newblock {\em J. Sci. Comp.}, 83:6, 2020.

\bibitem{Shu:97}
C.-W. Shu.
\newblock Essentially {N}on-{O}scillatory and {W}eighted {E}ssentially
  {N}on-{O}scillatory {S}chemes for {H}yperbolic {C}onservation {L}aws.
\newblock {\em NASA/CR-97-206253 ICASE Report No.97-65}, 1997.

\bibitem{2017Tavelli_SemiImplicitAllMach}
M.~Tavelli and M.~Dumbser.
\newblock A pressure-based semi-implicit space-time discontinuous {G}alerkin
  method on staggered unstructured meshes for the solution of the compressible
  {N}avier-{S}tokes equations at all {M}ach numbers.
\newblock {\em J. Comput. Phys.}, 341:341--376, 2017.

\bibitem{VTSP23:Quinpi:Book}
G.~Visconti, S.~Tozza, M.~Semplice, and G.~Puppo.
\newblock A conservative a-posteriori time-limiting procedure in {Q}uinpi
  schemes.
\newblock In G.~Albi, W.~Boscheri, and M.~Zanella, editors, {\em Advances in
  Numerical Methods for Hyperbolic Balance Laws and Related Problems. {YR}
  2021}, volume~32 of {\em SEMA-SIMAI Springer Series}, pages 191--212.
  Springer Cham, 2023.

\bibitem{ZDLS:14}
O.~Zanotti, M.~Dumbser, R.~Loub\`{e}re, and S.~Diot.
\newblock A posteriori subcell limiting for discontinuous {G}alerkin finite
  element method for hyperbolic system of conservation laws.
\newblock {\em J. Comput. Phys.}, 278:47--75, 2014.

\end{thebibliography}

\end{document}